\documentclass[12pt]{article}
\usepackage[a4paper]{geometry}
\usepackage{amsthm}
\usepackage{amsmath}
\usepackage{amssymb}
\usepackage{amsfonts}
\usepackage{graphicx}
\numberwithin{equation}{section}
\def\eref#1{(\ref{#1})}
\def\N{\mathbb{N}}
\def\P{\mathbb{P}}
\def\R{\mathbb{R}}
\def\T{\mathbb{T}}
\def\S{\mathbb{S}}
\def\E{\mathbb{E}}

\def\Z{\mathbb{Z}}

\def\C{\mathcal{C}}
\def\Pr{\mathcal{P}}

\def\F{\mathcal{F}}

\def\<{\big\langle}
\def\>{\big\rangle}
\def\Osc{\operatorname{Osc}}

\def\dist{\operatorname{dist}}

\newtheorem{Lemma}{Lemma}[section]
\newtheorem{Corollary}{Corollary}[section]
\newtheorem{Theorem}{Theorem}[section]
\newtheorem{Proposition}{Proposition}[section]
\newtheorem{Conjecture}{Conjecture}[section]

\newtheorem{Condition}{Condition}[section]
\theoremstyle{remark}
\newtheorem{Remark}{Remark}[section]
\newtheorem{Example}{Example}[section]
\theoremstyle{definition}

\newtheorem{Definition}{Definition}[section]

\begin{document}

\title{Multi-scale homogenization with bounded ratios and Anomalous 
Slow Diffusion.\protect\footnotetext{AMS 1991 {\it{Subject 
Classification}}. Primary 74Q20, 60J60; secondary   35B27, 74Q10,  60F05, 31C05, 
35B05.} \protect}          % Enter your title between 

\author{G\'{e}rard Ben Arous\footnote{ gerard.benarous@epfl.ch, DMA,  
EPFL, CH-1015 Lausanne , Switzerland} and Houman Owhadi\footnote{
 houman.owhadi@epfl.ch, DMA,  EPFL, CH-1015 Lausanne , Switzerland}}        
% Enter your name between curly braces

      % Enter your date or \today between curly braces
\maketitle
\abstract{We show that the effective diffusivity matrix $D(V^n)$ for the heat operator $\partial_t-(\Delta/2-\nabla V^n \nabla)$ in a periodic potential  $V^n=\sum_{k=0}^n 
U_k(x/R_k)$ obtained as a superposition of Holder-continuous periodic potentials  $U_k$ (of period 
$\T^d:=\R^d/\Z^d$, $d\in \N^*$, $U_k(0)=0$) decays exponentially fast with the number of scales when the scale-ratios $R_{k+1}/R_k$ are bounded above and below. From this we deduce 
the anomalous slow behavior for a Brownian Motion in a potential obtained as a superposition of an infinite number of scales: $dy_t=d\omega_t -\nabla V^\infty(y_t) dt$}

\tableofcontents

\section{Introduction}
Homogenization in the presence of a large number of spatial scales is 
both very important for applications and far from understood from a 
mathematical standpoint. In the asymptotic regime where the spatial scales 
separate, i.e. when the ratio between successive scales tends to 
infinity, multi-scale homogenization is now well understood. See for instance 
Bensoussan-Lions-Papanicolaou \cite{BeLiPa78}, or Avellaneda 
\cite{Av87}, Allaire - Briane \cite{AlBr96}, \cite{Ko95}, Jikov - Kozlov 
\cite{JiKo99} or Avellaneda - Majda \cite{AvMa94}.\\
Nevertheless the case of multi-scale homogenization when spatial scales 
are not clearly separated, i.e. when the ratios between scales stay 
bounded, has been recognized as difficult and important. For instance, 
Avellaneda \cite{Ave96} (page 267) emphasizes that "the assumption of 
scale separation invoked in homogenization is not adequate for treating the 
most general problems of transport and diffusion in self-similar random 
media".\\
The potential use of multi-scale homogenization estimates for 
applications are numerous (see for instance \cite{SeScCo81} for applications to 
geology, or \cite{Bou76}, \cite{McL77}, \cite{ClChLe80} for 
applications to Differential Effective Medium Theories).
The main application of this line of ideas is perhaps to proving 
super-diffusivity for turbulent diffusion: see for instance \cite{Ave96},  
\cite{AvMa90};  \cite{FuGl90}, \cite{FuGl91}, \cite{Gli92}, 
\cite{GlZh92}, \cite{Zh92},  \cite{IsKa91},  \cite{FaPa94},\cite{FaPa96}, 
\cite{Fan99};  \cite{Bha99}, \cite{BhDeGo99}; \cite{FaKo01} or \cite{KoOl00}.\\
We are here interested in sub-diffusivity problems. Consider the 
Brownian motion in a periodic potential, i.e. the diffusion process
\begin{equation}
dy_t=d\omega_t-\nabla V(y_t) dt
\end{equation}
where $V$ is periodic and smooth. It is a basic and simple fact of 
homogenization theory that $y_t$ behaves in large times like a Brownian 
motion slower than the Brownian Motion $\omega_t$ driving the equation, 
i.e. $y^\epsilon(t)=\epsilon y_{t/\epsilon^2}$ converges in law to a 
Brownian motion with diffusivity matrix $D(V)<I_d$. \\
We first treat here the case where $V$ is a periodic $n$-scale 
potential with ratios (between successive scales) bounded uniformly on $n$.We 
introduce  a new approach which enables us to show exponential decay of 
the effective diffusivity matrix when the number of spatial scales 
grows to infinity.\\
From this exponential decay we will deduce the anomalous slow behavior 
of Brownian motions in potential $V$, when $V$ is a superposition of an 
infinite number of scales.\\
We have studied this question with a particular application in mind, 
i.e. to prove that one of the basic mechanisms of anomalous slow 
diffusion in complex media is the existence of a large number of spatial 
scales, without a clear separation between them.  This phenomenon has been 
attested for very regular self-similar fractals (see Barlow and Bass 
\cite{BB97} and Osada \cite{Os90} for the Sierpinski carpet) (see also  
\cite{HaKu98}). Our goal is to implement rigorously  the idea that the key 
for the sub-diffusivity is a never-ending or perpetual homogenization 
phenomenon over an infinite number of scales, the point being that our 
model will not have any self similarity or local symmetry hypotheses.\\
Our approach gives naturally much more detailed information in 
dimension one and this is the subject of \cite{Ow00a}.\\ This approach will be 
shown in forthcoming works to also give a proof of super-diffusive 
behavior for diffusion in some multi-scale divergence free fields  (see 
\cite{BeOw00c} for the simple case of shear flow and \cite{Ow00b} for a 
general situation).\\
The second section contains the description of our model; the third 
one, the statement of our results and the fourth one the proofs.

\nocite{MajKra99} \nocite{Ol94}

\section{The multi-scale medium}
\begin{figure}
\includegraphics[scale=0.8]{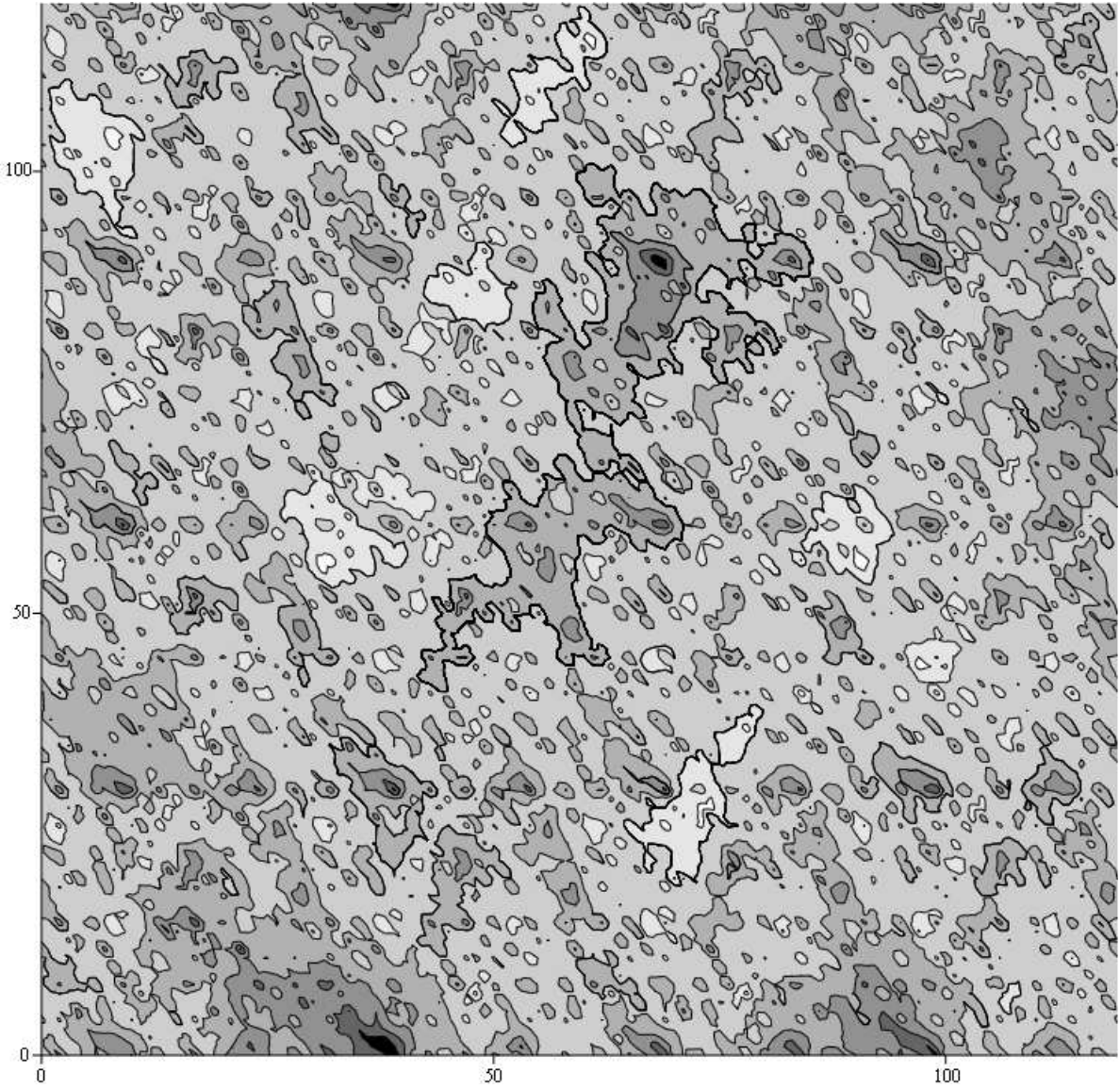}
\caption{A particular case \label{Figspf1}}
\end{figure}

 For $U \in L^\infty(\T^d_R)$ (we note $\T^d_R := R \T^d$), let $m_U$ 
be the probability measure on $\T^d$ defined by
\begin{equation}
m_U(dx)=e^{-2U(x)}dx\big/\int_{\T^d_R}e^{-2U(x)}dx
\end{equation}
The effective diffusivity $D(U)$ is the symmetric positive definite 
matrix given by
\begin{equation}\label{eqforduvarfjhoiuh71}
^tlD(U)l=\inf_{f\in C^\infty(\T^d_R)}\int_{\T^d_R} |l-\nabla f(x)|^2 
m_U(dx)
\end{equation}
for $l$ in $\S^{d-1}$ (the unit sphere of $\R^d$).
Our purpose in this work is to obtain quantitative estimates for the 
effective diffusivity matrix of multi-scale potentials $V_0^n$  given by 
a sum of  periodic functions with
(geometrically) increasing periods:
\begin{equation}\label{Modsubfracuinfty}
V_0^n=\sum_{k=0}^n U_k(\frac{x}{R_k})
\end{equation}
In this formula we have two important ingredients: the potentials $U_k$ 
and the scale parameters $R_k$. We will now describe the hypothesis we 
make on these two items of our model.
\begin{enumerate}
\item \underline{Hypotheses on the potentials $U_k$}\\
We will assume that  
\begin{equation}\label{jhshdddikuou1}
U_k \in C^\alpha(\T^d)
\end{equation}
\begin{equation}\label{jhshdikuou1}
U_k(0)=0
\end{equation}
Here $C^\alpha(\T^d)$ denotes the space of $\alpha$-Holder continuous  
on the torus $\T^d$, with $0<\alpha\leq 1$. 
We will also assume that the $C^\alpha$-norm of the $U_k$ are uniformly 
bounded, i.e.
\begin{equation}\label{ModsubContUngradUn}
K_\alpha:=\sup_{k\in \N} \sup_{x\not=y}|U_k(x)-U_k(y)|/|x-y|^\alpha 
<\infty
\end{equation}
We will also need the notation
\begin{equation}\label{ModsubContUngradUnGer}
K_0:=\sup_{k\in \N} \Osc(U_k) 
\end{equation}
where the oscillation of $U_k$ is given by $\Osc(U):=\sup U - \inf 
U$.\\
We also assume that the effective diffusivity matrices of the $U_k$'s 
are uniformly bounded. Let $\lambda_{\min}(D(U_k))$ (respectively 
$\lambda_{\max}(D(U_k))$) be the smallest and largest eigenvalues of the 
effective diffusivity matrix $D(U_k)$. We will assume that
\begin{equation}\label{ModsubDiffCondUniDUn}
\lambda_{\min}:=\inf_{k\in \N, l\in S^d} {^tlD(U_k)l}>0
\end{equation}
\begin{equation}\label{ModsubDiffCondUniDUnGer}
 \lambda_{\max}:=\sup_{k\in \N, l\in S^d} {^tlD(U_k)l}<1
\end{equation}
\item \underline{Hypotheses on the scale parameters $R_k$}\\
$R_k$ is a spatial scale parameter growing exponentially fast with $k$, 
more precisely we will assume that $R_0=r_0=1$ and that the ratios 
between scales defined by
\begin{equation}\label{jahgsvagvjh6761}
 r_k= R_k/R_{k-1}\in \N^{*}
\end{equation}
for $k\geq 1$, are integers uniformly bounded away from $1$ and 
$\infty$: we will denote by
\begin{equation}\label{Modsubboundrnrhonmin}
 \rho_{\min}:=\inf_{k\in \N^*} r_k \quad\text{and}\quad 
\rho_{\max}:=\sup_{k\in \N^*} r_k 
\end{equation}
and assume that 
\begin{equation}\label{Modsubboud32ndrnrhonmin}
 \rho_{\min}\geq 2 \quad\text{and}\quad \rho_{\max}<\infty 
\end{equation}
\end{enumerate}
As an example, we have illustrated in the figure \ref{Figspf1}  the 
contour lines of
$V_0^2(x,y)=\sum_{k=0}^2 U(\frac{x}{\rho^k},\frac{y}{\rho^k})$
, with $\rho=4$ and
$U(x,y)=\cos(x+\pi\sin(y)+1)^2\sin(\pi \cos(x)-2y+2)\cos(\pi\sin 
(x)+y)$

\section{Main results}
\subsection{Quantitative estimates of the multi-scale effective 
diffusivity}\label{wsqlkiuuw771sd}
\subsubsection{The central estimate}
Our first objective is to control the minimal and maximal eigenvalues 
of $D(V_0^n)$. More precisely writing $I_d$ the $d\times d$ identity 
matrix we will prove that
\begin{Theorem}\label{Masubthentr37284gd1}
Under the hypotheses \eref{ModsubContUngradUn}, \eref{jahgsvagvjh6761} 
and $\rho_{\min}^\alpha \geq K_\alpha$  there exists a constant $C$ 
depending only on $d,\alpha,K_\alpha,K_0$ such that for all $n \geq 1$
\begin{equation}\label{Thupboprthuv81}
I_d\; e^{-n \epsilon}\prod_{k=0}^{n}\lambda_{\min}\big(D(U_k)\big)
\leq D(V_0^{n}) \leq I_d\; e^{n \epsilon} \prod_{k=0}^{n} 
\lambda_{\max}(D(U_k))
\end{equation}
where 
\begin{equation}\label{suhziszsz771}
\epsilon=C \rho_{\min}^{-\alpha/2}
\end{equation}
In particular $\epsilon$ tends to $0$, when $\rho_{\min}\rightarrow 
\infty$
\end{Theorem}
\begin{Remark}
One can interpret this theorem as follows: $D(V_0^{n})$ is bounded from 
below and from above by  the bounds given by reiterated homogenization 
under the assumption of complete separation between scales, i.e. 
$\rho_{\min}\rightarrow \infty$ (product of minimal and
maximal eigenvalues) times an error term $e^{n \epsilon}$ created by 
the \emph{interaction} or \emph{overlap} between the different scales.
\end{Remark}
\begin{Remark}
Originally the problem of estimating $D(V_0^n)$ was called for in 
connection with applied sciences,  and heuristic theories such as 
Differential Effective Medium theory  have been developed for that purpose. This 
theory (DEM theory) models a two phase composite by incrementally 
adding inclusions of one phase to a background matrix of the other and then 
recomputing the new effective background material at each increment 
\cite{Bou76}, \cite{McL77}, \cite{ClChLe80}. It was first proposed by 
Bruggeman to compute the conductivity of a two-component composite 
structure formed by successive substitutions (\cite{Bru35} and \cite{AIP77}) 
and generalized by Norris \cite{Nor85} to materials with more than two 
phases.\\
More recently Avellaneda \cite{Av87} has given a rigorous 
interpretation of the equations obtained by DEM theories showing that they are 
\emph{homogeneous limit equations} with two very important features: 
complete separation of scales and "dilution of phases". That is to say, each 
"phase" $U_k$ is present at an infinite number of scales in a 
homogeneous way. Yet two different phases never interact because they always 
appear at scales whose ratio is $\infty$. Moreover the macroscopic 
influence of each phase is totally (but non uniformly) diluted in the infinite 
number of scales at which it appears. In our context, complete 
separation of scales would mean that $R_{k+1}/R_k$ grows sufficiently fast to 
$\infty$, and "dilution of phases" would bmean that $V_0^n=\sum_{k=0}^n 
U_k^n(x/R_k)$ with $U_k^n\rightarrow 0$ as $n\rightarrow \infty$. The 
rigorous tool used by Avellaneda to obtain this interpretation is 
reiterated homogenization  \cite{BeLiPa78}.\\ 
May be the most recent work on this topic is the article by 
Jikov-Kozlov \cite{JiKo99}, who work under the assumption of "dilution of phases" 
and  fast separation between scales, more precisely under the condition 
that $\sum_{k=1}^\infty k\big(R_k/R_{k+1}\big)^2<\infty$.
Jikov-Kozlov use the classical toolbox of asymptotic expansion, 
plugging well chosen test functions in the cell problem. This method of 
asymptotic expansion is simply not at all available in our context.
\end{Remark}
\bigskip
Theorem \ref{Masubthentr37284gd1} will be proved by induction on the 
number of scales.
The basic step in this induction is the following estimate 
\eref{Thupboprthuv82} on the effective diffusivity for a two-scale periodic 
medium.\\
Let $U, T\in C^\alpha(\T^d)$. Let us define for $R\in \N^*$, $S_R U \in 
C^\alpha(\T^d)$ by $S_R U(x)=U(Rx)$. We will need to estimate $D(S_R 
U+T)$ the effective diffusivity for a two-scale medium when $R$ is a 
large integer.
Let us define $D(U,T)$, the symmetric definite positive matrix given by
\begin{equation}\label{eqforduvarfjhoiuh71jGer}
^tlD(U,T)l=\inf_{f\in C^\infty(\T^d_1)}\int_{\T^d_1} {^t(}l-\nabla 
f(x))D(U)(l-\nabla f(x)) m_T(dx)\quad \text{for $l\in \R^d$}
\end{equation}
\begin{Theorem}\label{NewthcontmscM}
Let $R\in \N^*$ and $U,T \in C^\alpha(\T^d)$. If $R^\alpha \geq 
\|T\|_\alpha$ then there exists a constant $C$ depending only on 
$d,\Osc(U),\|U\|_\alpha,\alpha$ such that
\begin{equation}\label{Thupboprthuv82}
e^{- \epsilon}
D(U,T) \leq D(S_R U+T)\leq D(U,T) e^{ \epsilon}
\end{equation}
with $\epsilon=C R^{-\alpha/2}$
\end{Theorem}
\begin{Remark}
 Theorem \ref{NewthcontmscM} implies obviously that
\begin{equation}
D(U,T)=\lim_{R\rightarrow \infty}D(S_R U+T)
\end{equation}
so that $D(U,T)$ should be interpreted as the effective diffusivity of 
the two scale medium for a complete separation of scales. Naturally 
$D(U,T)$ is also computable from an explicit cell problem (see 
\eref{sjddskdjhv675A}).
\end{Remark}
\begin{Remark}
 The estimate given in theorem \ref{NewthcontmscM} is stronger than 
needed for theorem \ref{Masubthentr37284gd1}. It gives a control of $D(S_R 
U+T)$ in terms of $D(U,T)$ and not only of the minimal and maximal 
eigenvalues of $D(U)$ and $D(T)$. In fact we will only use its corollary 
\ref{Nej88871wthcontmscM} given below, which is deduced using the 
variational formulation \eref{eqforduvarfjhoiuh71jGer}.
\end{Remark}
\begin{Corollary}\label{Nej88871wthcontmscM}
Let $R\in \N^*$ and $U,T \in C^\alpha(\T^d)$. If $R^\alpha \geq 
\|T\|_\alpha$ then there exists a constant $C$ depending only on 
$d,\Osc(U),\|U\|_\alpha,\alpha$ such that
\begin{equation}\label{Thupboprthuv82Ger}
\lambda_{\min}\big(D(U)\big)
D(T)e^{- \epsilon} \leq D(S_R U+T)\leq \lambda_{\max}\big(D(U)\big)D(T) 
e^{ \epsilon}
\end{equation}
with $\epsilon=C R^{-\alpha/2}$
\end{Corollary}
\begin{Remark}
We mentioned that theorem \ref{Masubthentr37284gd1} is proved by 
induction. This induction differs from the one used in reiterated 
homogenization or DEM theories by the fact we homogenize on the larger scales 
first and add at each step a smaller scale.
\end{Remark}
Let us introduce the following upper and lower exponential rates
\begin{Definition}
\begin{equation}
\lambda^+=\lim\sup_{n\rightarrow \infty}\frac{1}{n} \ln 
\lambda_{\max}\big(D(V_0^n)\big)
\end{equation}
\begin{equation}
\lambda^-=\lim\inf_{n\rightarrow \infty}\frac{1}{n} \ln 
\lambda_{\min}\big(D(V_0^n)\big)
\end{equation}
\end{Definition}
Theorem \ref{Masubthentr37284gd1} implies the exponential decay of 
$D(V_0^{n-1})$, i.e.
\begin{Corollary}\label{Masubthentr37284gd2}
Under the hypotheses \eref{ModsubContUngradUn}, \eref{jahgsvagvjh6761} 
and $\rho_{\min}^\alpha\geq K_\alpha$ one has (with $\epsilon$ given by 
\eref{suhziszsz771}) for $n\geq 1$
\begin{equation}\label{eqesdvo0n61}
I_d e^{-n \epsilon}\lambda_{\min}^{n+1}
\leq D(V_0^{n}) \leq I_d e^{n \epsilon} \lambda_{\max}^{n+1}
\end{equation}
and
\begin{equation}
\lambda^+\leq \ln \lambda_{\max}+\epsilon
\end{equation}
\begin{equation}
\lambda^-\geq \ln \lambda_{\min}-\epsilon
\end{equation}
In particular if $\lambda_{\max}<1$ then there exists a constant 
$\rho_0=\big(1+C_{d,K_0,K_\alpha,\alpha}/(-\ln 
\lambda_{\max}))^\frac{2}{\alpha}$ such that, for $\rho_{\min} \geq \rho_0$
\begin{equation}
\lambda^+<0
\end{equation}
\end{Corollary}
Thus one obtains  the exponential decay of $D(V^n)$  only for a minimal 
separation between scales, i.e. $\rho_{min}$ greater than a constant 
$\rho_0$ characterized by the medium. It is natural to wonder whether 
this condition is necessary  and what happens below this constant 
$\rho_0$. We will give a partial answer to that question, in the simple case 
when the medium $V$ is self-similar. We will see that it is possible to 
find models such that for a certain value $C$ of the separation 
parameter $\rho_{\min}=C$, $D(V_0^n)$ decays exponentially and for 
$\rho_{\min}=C+1$, $D(V_0^n)$ stays bounded away from zero. This will be done using 
a link with large deviation theory.

\subsubsection{The self similar case}
\begin{Definition}
The medium $V$ is called self similar if and only if $\forall n$, 
$U_n=U$ and $R_n=\rho^n$ with $\rho\in \N$, $\rho\geq 2$
\end{Definition}
\begin{Definition}
For $U\in C^\alpha(\T^d)$ and $\rho\in \N/\{0,1\}$ we denote by 
$p_\rho(U)$ the  pressure associated to the shift $s_\rho(x)=\rho x$ on 
$\T^d$, i.e.
\begin{equation}
p_\rho(U)=\sup_{\mu}\Big(\int_{\T^d} U(x) d\mu(x)+h_\rho(\mu) \Big)
\end{equation}
where $h_\rho$ is the Komogorov-Sinai entropy related to the shift 
$s_\rho$. We denote
\begin{equation}
P_\rho(U)=p_\rho(U)-p_\rho(0)=p_\rho(U)-d\ln R
\end{equation}
\end{Definition}
We refer to \cite{Ke98} and \cite{Rue78} for a reminder on the 
pressure, let us observe that $\Pr_\rho(0)$ differs from the standard 
definition of the topological pressure by the constant $d \ln \rho$ so that 
$\Pr_\rho (0)=0$.\\
We will relate in the self-similar case the exponential rates 
$\lambda^+$ and $\lambda^-$ to  pressures for the shift $s_\rho$ and to large 
deviation at level $3$ for i.i.d. random variables.\\
In the self similar case we will write $\lambda^-(-U)$ the exponential 
rates associated to $D(-V_0^n)$. We write
\begin{equation}
Z(U)=-\Big(\Pr_\rho(2U)+\Pr_\rho(-2U)\Big)
\end{equation}

\begin{Theorem}\label{assjhsao78982hh1}
If the medium $V$ is self similar then
\begin{enumerate}
\item If $d=1$
\begin{equation}\label{eqjhl9889b87k1}
\lambda^+(U)=\lambda^-(U)=Z(U)
\end{equation}
\item\label{assssskavuuszz212Gert1} if $d=2$  then
\begin{equation}\label{assssskavuuszz}
\lambda^+(U)+\lambda^-(-U)=Z(U)
\end{equation}
Moreover if there exists an isometry $A$ of $\R^d$ such that 
$U(Ax)=-U(x)$ and a reflection $B$ such that $U(Bx)=U(x)$ then
$\lambda^-(-U)=\lambda^-(U)=\lambda^+(U)$ so that
\begin{equation}\label{assssskavuuszz212}
\lambda^+(U)=\lambda^-(U)=Z(U)/2
\end{equation}
\item For any $d$
\begin{equation}\label{dssdddjsdj5365}
Z(U)\leq \lambda^-(U) \leq \lambda^+(U)\leq 0
\end{equation}
\end{enumerate}
\end{Theorem}
\begin{Remark}
The statement \eref{eqjhl9889b87k1} is obtained from the explicit 
formula for $D(V_0^n)$ in $d=1$, see \cite{Ow00a}. 
\end{Remark}
\begin{Remark}
It is obviously important to know when $Z(U)$ is strictly negative to 
be able to use this theorem. A well known and useful criterion
can be stated as $Z(U)<0$ if and only if $U$ does not belong to the 
closure of the vector space spanned by cocycles, which can be shown to be 
equivalent to say that
\begin{equation}\label{eqjhl9889b87k2}
Z(U)<0 \Leftrightarrow \lim \sup_{n\rightarrow \infty} 
\frac{1}{n}\Big\|\sum_{k=0}^{n-1}\big(U(\rho^k 
x)-\int_{\T^d}U(x)dx\big)\Big\|_{\infty}>0
\end{equation}
We refer to \cite{Ow00a} for the proof of the last statement.
\end{Remark}
\begin{Example}
Let $U(x)=\sin(x)-\sin(81x)$ in dimension one. In fact 
\eref{eqjhl9889b87k1} and \eref{eqjhl9889b87k2} shows that $\lambda^+(U)<0$  as soon as 
$\rho \geq 82$. For $\rho\leq 81$ the situation is a bit surprising: 
$\lambda^+(U)<0$
for $\rho\not= 3,9,81$. For these exceptional values $\lambda^+(U)=0$ 
and in fact $D(V_0^n)$ remains lower bounded by a strictly positive 
constant.\\
This example shows that for a given potential $U$, even though the 
multi-scale effective diffusivity $D(V_0^n)$ decays exponentially for 
$\rho$ large enough, one can find isolated values of the scale parameter for 
which $D(V_0^n)$ remains bounded from below.   
\end{Example}
\begin{Remark}
The symmetry hypotheses given in theorem  
\ref{assjhsao78982hh1}.\ref{assssskavuuszz212Gert1} are only used to prove that for all $n$, 
$D(V_0^n)=D(-V_0^n)$ and 
$\lambda_{\max}\big(D(V_0^n)\big)=\lambda_{\min}\big(D(V_0^n)\big)$ (see proposition \ref{Thjbztvjhvzu7h71})
\end{Remark}

\subsection{Sub-diffusive behavior from homogenization on infinitely 
many scales}
Here we consider the diffusion process given by the Brownian Motion in 
the potential
\begin{equation}
V=V_0^\infty=\sum_{k=0}^\infty U_k(x/R_k)
\end{equation}
We assume in this section that the hypotheses \eref{jhshdddikuou1} , 
\eref{jhshdikuou1}, \eref{ModsubContUngradUn}, 
\eref{ModsubDiffCondUniDUn}, \eref{ModsubDiffCondUniDUnGer}, \eref{jahgsvagvjh6761} and 
\eref{Modsubboud32ndrnrhonmin} hold. To start with we will assume that 
\begin{equation}\label{sajhskz71}
\alpha=1 \quad \text{and that the potentials $U_k$ are uniformly 
$C^1$.} 
\end{equation}
In particular $V$ is well defined and belongs to $C^1(\R^d)$ and 
$\|\nabla V\|_\infty<\infty$.\\
The diffusion process associated to the potential $V$ is well defined 
by the Stochastic Differential Equation
\begin{equation}\label{IntModelsubdiffstochdiffequ}
  dy_t = d\omega_t - \nabla V(y_t) dt
\end{equation}
We will show that the multi-scale structure of $V$ can lead to an 
anomalous slow behavior for the process $y_t$. To describe this 
sub-diffusive phenomenon we choose to compute the mean exit time from large balls, 
i.e.
Let
\begin{equation}
\tau(r)=\inf\{t>0\,:\,|y_t|\geq r\}
\end{equation}
We would like to show that $\E_x[\tau(r)]$ grows faster than quadratic 
in $r$ when $r\rightarrow \infty$ uniformly in $x$. We cannot obtain 
such pointwise results in dimension $d>1$ (see \ref{dsjhdbbz817} for a 
discussion, the case $d=1$ is treated in \cite{Ow00a}). But we will start 
with averaged results on those mean exit times.\\
The fact that the homogenization results of subsection 
\ref{wsqlkiuuw771sd} can be of some help to estimate the mean exit times is shown by 
the following lemma
\begin{Lemma}\label{ksajhhdzz8717}
For $U\in C^\infty(\T^d_R)$, ($R>0$) writing $\E^{U}$ the exit times 
associated to the diffusion generated by $L_U=\Delta/2-\nabla U\nabla$ 
one has
\begin{equation}\label{Aldtau_alphaeq3}
\begin{split}
\E_x^U[\tau(x,r)]& \leq C_2 
\frac{r^2}{\lambda_{\max}\big(D(U)\big)}+C_d e^{(9d+15)\Osc(U)} R^2\\&\geq
C_1 \frac{r^2}{\lambda_{\max}\big(D(U)\big)}-C_d e^{(9d+15)\Osc(U)} R^2
\end{split}
\end{equation}
\end{Lemma}
\noindent Let $m_{V,r}$ be the probability measure on the ball $B(0,r)$ 
given by
\begin{equation}
m_{V,r}(dx)=e^{-2V(x)}\,dx\Big/\big(\int_{B(0,r)}e^{-2V(x)}\,dx\big)
\end{equation}
We will consider the mean exit time for the process started with 
initial distribution $m_{V,r}$, i.e.
\begin{equation}\label{aljhbdebii8918}
\E_{m_{V,r}}\big[\tau(r)\big]=\int_{B(0,r)}\E_x\big[\tau(r)\big]\,m_{V,r}(dx)
\end{equation}
\begin{Theorem}\label{IntSMAldsubanhianmethboaj1}
Under the hypothesis $\lambda_{\max}<1$ there exists $C_2$ depending on 
$d,\lambda_{\max},K_0,K_\alpha,\alpha$
such that if $\rho_{\min}>C_2$ then
\begin{equation}
\lim\inf_{r\rightarrow \infty} \frac{\ln 
\E_{m_{V,r}}\big[\tau(r)\big]}{\ln r}>2
\end{equation}
More precisely there exists $C_3>0$, $C_4>0$, $C_5>0$ such that for 
$r>C_{3}$,
\begin{equation}\label{hgvscdklvxkj1}
\E_{m_{V,r}}\big[\tau(r)\big]=r^{2+\nu(r)}
\end{equation}
with
\begin{equation}
0<C_4<\frac{\ln \frac{1}{\lambda_{\max}}}{\ln \rho_{\max}} 
\big(1-\frac{C_{5}}{\ln \rho_{\min}}\big)-\frac{1}{\ln r}C_{5} \leq \nu(r)\leq 
\frac{\ln \frac{1}{\lambda_{\min}}}{\ln \rho_{\min}} \big(1+\frac{C_5}{\ln 
\rho_{\min}}\big)
+\frac{1}{\ln r}C_5
\end{equation}
Where $C_3$ and $C_5$ depend on $(d,K_0,K_\alpha,\alpha)$ and $C_4$ on 
$(\lambda_{\max},\rho_{\max})$.
\end{Theorem}
The proof of this result relies heavily on theorem 
\ref{Masubthentr37284gd1}. The idea being that $\E_{m_{V,r}}\big[\tau(r)\big]$ is close, 
when $r$ is large to $r^2/ \lambda_{\max}\big(D(V_0^n)\big)$ where $n$ is 
roughly $\sup\{m\in \N\,:\, R_m\leq r\}$. So that the exponential decay 
of $D(V_0^n)$ gives the super-quadratic behavior of 
$\E_{m_{V,r}}\big[\tau(r)\big]$, i.e. sub-diffusivity.
\begin{Remark}
The differentiability hypothesis \eref{sajhskz71} though convenient in 
order to define the process $y_t$ as a solution of the SDE 
\ref{IntModelsubdiffstochdiffequ} is in fact useless. The theorem is also 
meaningful and true with $0<\alpha<1$.  See section \ref{kvdiudd22} for an 
explanation.
\end{Remark}

\subsubsection{Pointwise estimates on the anomaly}\label{dsjhdbbz817}
Theorem \ref{IntSMAldsubanhianmethboaj1} gives the anomalous behavior 
of the exit times with respect to the invariant measure of the diffusion 
and it is desirable to seek for pointwise estimates of this anomaly. 
The additional difficulty is to obtain quantitative estimates on the 
stability of divergence form elliptic operators under a perturbation of 
their principal parts (see conjecture \ref{sidusdviUZDSU}). By stability 
we mean here the validity of the following condition 
\ref{IntSubModCoStCojh1}.\\
For  $U \in C^1(\overline{B(z,r)})$. Write $\E^{U}$,  the expectation 
associated to the diffusions generated by 
$L_{U}=\frac{1}{2}\Delta-\nabla U\nabla$. $V$ is said to satisfy the stability condition 
\ref{IntSubModCoStCojh1} if and only if:
\begin{Condition}\label{IntSubModCoStCojh1}
There exists $\mu>0$ such that for all $n \in N$, all $z \in \R^d$, and 
all $r>0$,
\begin{equation}
\frac{1}{\mu} e^{-\mu \Osc_{B(z,r)} (V_{n+1}^\infty)} \inf_{x\in 
B(z,\frac{r}{2})} E_x^{ V_0^n} \big[\tau(B(z,r))\big]
 \leq E_z^{V} \big[\tau(B(z,r))\big]
\end{equation}
\begin{equation}
 E_z^{V} \big[\tau(B(z,r))\big] \leq \mu e^{\mu \Osc_{B(z,r)} 
(V_{n+1}^\infty)} \sup_{x\in B(z,r)} E_x^{ V_0^n} \big[\tau(B(z,r))\big]
\end{equation}
\end{Condition}
Where $\Osc_{B(z,r)} (U)$ stands for $\sup_{B(z,r)} U - \inf_{B(z,r)} 
U$.
Under the stability condition \ref{IntSubModCoStCojh1}, we can obtain 
sharp pointwize estimates on the mean exit times. 
\begin{Theorem} \label{IntSMAldcontrol_taudx2}
If $V$ satisfies the stability condition \ref{IntSubModCoStCojh1}, then 
there exist a constant $C_6$ depending on 
$(d,K_0,K_\alpha,\alpha,\mu,\lambda_{\max})$ such that for $\rho_{\min}>C_{6}$, one has for all $x 
\in \R^d$
\begin{equation}
\lim \inf_{r\rightarrow \infty} \frac{\ln \E_x 
\big[\tau(B(x,r))\big]}{\ln r}>2
\end{equation}
More precisely there exists a function $\sigma(r)$ such that for
$r>C_{7}$ one has
\begin{equation}
  C_{8} r^{2+\sigma(r)(1-\gamma)} \leq \E_x \big[\tau(B(x,r))\big] \leq 
C_{9}   r^{2+\sigma(r)(1+\gamma)}
\end{equation}
with
\begin{equation}
\frac{\ln \frac{1}{ \lambda_{\max}}}{\ln \rho_{\max}} 
(1+\frac{C_{3}}{\ln \rho_{\min}})^{-1}
 \leq \sigma(r)\leq \frac{\ln \frac{1}{ \lambda_{\min}}}{\ln 
\rho_{\min}} (1+\frac{C_{4}}{\ln \rho_{\min}})
\end{equation}
and $\gamma=C_{5}K_0\big/(\ln \rho_{\min})<0.5$. Where the constants 
$C_3,C_4,C_7,C_8$ and $C_9$ depend on 
$(d,K_0,K_\alpha,\alpha,\mu,\lambda_{\max})$ and $C_5$ on $d$.
\end{Theorem}
\noindent Here $\tau(B(x,r))$ denotes the exit time from the ball 
$B(x,r)$.
\begin{Remark}
In fact $\sigma(r)$ can be described rather precisely. Let,
\begin{equation}
\sigma(r,n)=-\ln \lambda_{\max} D(V_0^n)\Big/\ln r
\end{equation}
Define $n_{ef}(r,C_1,C_2)=\sup\{n \geq 0 \;:\;  e^{(n+1)C_{1}K_0} R^2_n  
\leq C_{2} r^2\}$. Then there exists $C_1,C_2$ depending only on $d$ 
such that $\sigma(r)$ in theorem \ref{IntSMAldcontrol_taudx2} is 
$\sigma\big(r,n_{ef}(r,C_1,C_2)\big)$.
\end{Remark}
Using the precise information of theorem \ref{IntSMAldcontrol_taudx2} 
we can estimate the tails of probability transitions for the process 
$y_t$ (or the tail of the heat kernel for the operator $L_V$). We get 
non-Gaussian upper bound similar to the (more precise) ones proved for 
fractal diffusions (see \cite{BB97} and \cite{HaKu98})

\begin{Theorem}\label{ddshdidcsbiii1}
If $V$ satisfies the stability condition \ref{IntSubModCoStCojh1}, then 
for $\rho_{\min}>C_{6}$, $r>0$
\begin{equation}\label{gvisuziuzvvvv1}
C_{10} r \leq t \leq C_{11} r^{2+\sigma(r)(1-3\gamma)}
\end{equation}
one has
\begin{equation}\label{skassljs771}
\ln \P_x[|y_t-x|\geq r]  \leq -C_{13} \frac{r^2}{t} 
\big(\frac{t}{r}\big)^{\nu(t/h)}
\end{equation}
with ($C_{17}<0.5 \ln \rho_{\min}$)
\begin{equation}
0<\frac{\ln \frac{1}{\lambda_{\max}}}{\ln 
\rho_{\max}}(1-\frac{C_{14}}{\ln \rho_{\min}}) \leq \nu(y) \leq \frac{\ln 
\frac{1}{\lambda_{\min}}}{\ln \rho_{\min}}(1-\frac{C_{15}}{\ln \rho_{\min}})
\end{equation}
Where $C_6,C_{13},C_{14}$ and $C_{15}$ depend on 
$(d,K_0,K_\alpha,\alpha,\mu,\lambda_{\max})$; $C_{17}$ on $d,K_0$; $C_{10},C_{11}$ on
$(d,K_0,K_\alpha,\alpha,\mu,\lambda_{\max},\rho_{\min},\rho_{\max})$
\end{Theorem}
\begin{Remark}
 The non Gaussian structure of \eref{skassljs771} is similar to the one 
obtained for diffusion processes in fractals. Indeed
\begin{equation}\label{esdqart1sugtapsierptxy1}
-C \frac{h^2}{t} (\frac{t}{h})^{\nu}= - C 
\big(\frac{|x-y|^{d_w}}{t}\big)^\frac{1}{d_w-1}
\end{equation}
with $d_w\sim 2+\nu$.\\
Next, it has been shown in \cite{Ow00a} that for $U\in L^\infty( 
\T^d_R)$, 
$\ln p^U(x,y,t)$ is roughly $-{^t(}y-x)D^{-1}(U)(y-x)/t$ for $t > R 
|x-y|$ (homogenized behavior). Where  $p^U$ is the heat kernel associated 
to $L_U$. Next, writing $n_{ef}(t/h)=\sup_{n}\{R_n \leq t/h\}$ the 
number of scales that one can consider as homogenized in the estimation of 
the heat kernel tail one obtains from a heuristic computation (which 
can be made rigorous in dimension one, see \cite{Ow00a}) that for $C_{10} 
h \leq t \leq C_{11} h^{2+\mu}$, 
\begin{equation}\label{eqart1sugtapsierptxy1}
\ln \P(y_t\geq h)\leq -C \frac{h^2}{t \lambda^{n_{ef}(t/h)}}\sim -C 
\frac{h^2}{t} (\frac{t}{h})^{-\frac{\ln \lambda}{\ln \rho}}\sim - C 
\big(\frac{|x-y|^{d_w}}{t}\big)^\frac{1}{d_w-1}
\end{equation}
with $d_w\sim 2- \frac{\ln \lambda}{\ln \rho}$. The equation 
\ref{eqart1sugtapsierptxy1} suggests that the origin of the anomalous shape of 
the heat kernel for the reflected Brownian Motion on the Sierpinski 
carpet can be explained by a \emph{perpetual homogenization} phenomenon and 
the formula linking the number of effective scales and the ratio 
$t/h$.\\
The  condition $C_{10} h \leq t$ can be translated into "homogenization 
has started on at least the first scale" ($n_{ef}\geq 1$) and the 
second one $t \leq C_{11} h^{2+\mu}$ into "the heat kernel associated to 
$L_V$ is far from its diagonal regime"  (one can have $h^2/t<<1$ before 
reaching that regime, this is explained by the slow down of the 
diffusion).
\end{Remark}
The weak point of theorems \ref{IntSMAldcontrol_taudx2} and 
\ref{ddshdidcsbiii1} is naturally that
checking condition \ref{IntSubModCoStCojh1} seems difficult. But we 
believe that in fact this condition is always true (we refer to the 
chapter 13 of \cite{Owh01}), since this condition is a consequence of the 
following conjecture  (see \cite{Ow00a}, proposition 2.3). 
\begin{Conjecture}\label{sidusdviUZDSU}
There  exists a constant $C_{d}$
depending only on the dimension of the space  such that
for $\lambda \in C^\infty(\overline{B(0,1)})$ such that $\lambda>0$ on 
$\overline{B(0,1)}$ and $\phi,\psi \in C^2(\overline{B(0,1)})$ null on 
$\partial B(0,1)$ and both sub harmonic with respect to the operator 
$-\nabla(\lambda \nabla)$, one has
\begin{equation}
\int_{B(0,1)}\lambda(x) |\nabla \phi(x).\nabla \psi(x)|\,dx \leq C_{d} 
\int_{B(0,1)}\lambda(x) \nabla \phi(x).\nabla \psi(x)\,dx
\end{equation}
\end{Conjecture}
It is simple to see  \cite{Ow00a} that this conjecture is true  in 
dimension one with $C_{d}=3$. So proving conjecture 
\ref{sidusdviUZDSU} would give the pointwise estimates of theorem 
\ref{IntSMAldcontrol_taudx2} and the tail estimate for the heat kernel in 
theorem \ref{ddshdidcsbiii1}.
\begin{Remark}
Here we have assumed that the $U_k$ are uniformly $C_1$ but let us 
observe that since theorems \ref{IntSMAldcontrol_taudx2} and  
\ref{ddshdidcsbiii1} are robust in their dependence on $\alpha$ and $K_\alpha$ (one 
can choose $\alpha<1$). One can built a process, with the assumption 
that the $U_k$ are $\alpha$-Holder continuous, whose mean exit times and 
heat kernel tail verify the estimates given in  theorems 
\ref{IntSMAldcontrol_taudx2} and  \ref{ddshdidcsbiii1}. 
\end{Remark}

\section{Proofs}
\subsection{Multi-scale homogenization with bounded ratios}
\subsubsection{Global estimates of the multi-scale ef{f}ective 
dif{f}usivity: theorem \ref{Masubthentr37284gd1}}
The proof of theorem \ref{Masubthentr37284gd1} will follow from the 
Corollary \ref{Nej88871wthcontmscM} by a simple induction. Let $n\in 
N/\{0,1\}$, $p\in \N, 1\leq p\leq n$ and assume that
\begin{equation}\label{Thupboprthuv8hnv1}
I_d e^{-(n-p) 
\epsilon(\rho_{\min})}\prod_{k=p}^{n}\lambda_{\min}\big(D(U_k)\big)
\leq D(V_p^{n}) \leq I_d e^{(n-p) \epsilon(\rho_{\min})} 
\prod_{k=p}^{n} \lambda_{\max}(D(U_k))
\end{equation}
One pass from the quantitative control on $D(V_{p}^{n})$ to a control 
on $D(V_{p-1}^{n})$ by choosing $U(x)=U_{p-1}(x)$, 
$T(x)=V_{p}^{n}(R_{n}x)$ and $R=R_{n}/R_{p-1}$ in theorem \ref{NewthcontmscM} and observing 
that
$\|T\|_\alpha/R^\alpha \leq  (2^\alpha-1)^{-1} 
K_\alpha/\rho_{\min}^\alpha$. This proves the induction and henceforth the theorem.

\subsubsection{Quantitative multi-scale-homogenization: Upper bound in 
the  theorem \ref{NewthcontmscM}}
\paragraph{}\label{hcsajscshsgcvzsvc1}
We will use the notation introduced in theorem \ref{NewthcontmscM}. By 
the variational formula \eref{eqforduvarfjhoiuh71}, $D(U)$ is 
continuous with respect to $U$ in $L^\infty$-norm,  it is sufficient to prove 
theorem  \ref{NewthcontmscM}
assuming that $U$ and $T$ are smooth.\\
First let us prove that when homogenization takes place on two scales 
separated by a ratio $R$, the influence of a translation of the first 
one with respect to the second one  on the global effective diffusivity 
can easily be controlled, i.e. for $y\in T^d_1$, writing $\Theta_y$ the 
translation operator $T(x)\rightarrow \Theta_y T(x)=T(x+y)$
\begin{Lemma}\label{jhibuhzbzusbb1}
\begin{equation}\label{dsgutgveq17h7100}
e^{-4\frac{\| T\|_\alpha}{R^\alpha}} D(S_RU+ T) \leq D(S_R \Theta_y U+ 
T)\leq e^{4\frac{\| T\|_\alpha}{R^\alpha}}D(S_R U+ T)
\end{equation}
\end{Lemma}
\begin{proof}
The proof  follows by observing that $S_RU+\Theta_y T=\Theta_{[Ry]/R} 
(S_RU+ T)+
\Theta_{y}T-\Theta_{[Ry]/R}T$ where $[Ry]$ stands for the vector with 
the integral parts of $(yR)_i$ as coordinates.
Thus by the variational definition of the effective diffusivity
\begin{equation}
D(S_RU+\Theta_y T)\leq 
e^{4\|\Theta_{y}T-\Theta_{[Ry]/R}T\|_\infty}D(\Theta_{[Ry]/R} (S_RU+ T))
\end{equation}
and the equation \eref{dsgutgveq17h7100} follows by observing that the 
effective diffusivity is invariant under a translation of the medium: 
$D(\Theta_{[Ry]/R} (S_RU+ T))=D(S_RU+ T)$.
\end{proof}
Next we will obtain a quantitative control  on $\int_{y \in \T^d } 
D(S_R U+\Theta_y T) dy$ 
\begin{Lemma}\label{jhibuhzbzusbb1a}
For $R>\|T\|_\alpha$
\begin{equation}\label{dsgutgveq17h7100a}
\int_{y \in \T^d } D(S_R U+\Theta_y T) dy \leq e^{ 22 \frac{\| 
T\|_\alpha}{R^\alpha}} \big(1+C_d e^{C_d \Osc(U)} 
(\|T\|_\alpha/R^\alpha)^\frac{1}{2}\big)  D(U,T)
\end{equation}
\end{Lemma}
Let us observe that the combination of lemma \ref{dsgutgveq17h7100} 
with \ref{jhibuhzbzusbb1} gives
the upper bound \eref{Thupboprthuv82} in theorem \ref{NewthcontmscM}.\\
Write $\chi^U_l$ the solution of the cell problem associated to $U$. We 
remind that for $l\in \R^d$,  $L_U=1/2\Delta-\nabla U\nabla$, $L_U 
\chi_l=-l\nabla U$, $\chi^U_l(0)=0$ and 
\begin{equation}
^tlD(U)l=\int_{\T^d}|l-\nabla\chi_l|^2 
m_U(dx)=\int_{\T^d}{^t(}l-\nabla\chi_l).l m_U(dx)
\end{equation}
Write $\chi^{D(U),T}$ the $\T^d$ periodic solution of the following 
cell problem (which corresponds to a complete homogenization on the 
smaller scale): for $l\in \S^{d-1}$
\begin{equation}
\nabla \big(e^{-2T}D(U)(l-\nabla \chi^{D(U),T}_l) \big)=0
\end{equation}
Write for $y \in \T^d$, $x\rightarrow\chi(x,y)$ the solution of the 
cell problem associated to $S_R \Theta_y U+T$.\\
Let $l\in \S^{d-1}$, by the formula associating the effective 
diffusivity and the solution of the cell problem and using that $l-\nabla_x 
\chi_l(x,y)$ is harmonic with respect to $L_{S_R U+\Theta_y T}$, one 
obtains
\begin{equation}\label{dsgutgveq17h71Ger}
\begin{split}
\int_{y\in \T^d}{^tlD}(S_R \Theta_y U+T)ldy=\int_{ \T^d\times 
\T^d}(l-\nabla_x \chi_l(x,y)).l \,dx\,dy
\end{split}
\end{equation}
Writing the decomposition 
\begin{equation}
l=(I_d-\nabla \chi^{U}_.(Rx+y))(l-\nabla \chi^{D(U),T}_l(x))+\nabla 
\chi^{U}_.(Rx+y)(\nabla \chi^{D(U),T}_l(x)-l)m_{U(R.x+y)+T(.)}(dx)\,dy
\end{equation}
we get that
\begin{equation}\label{dsgutgveq17h71}
\begin{split}
\int_{y\in \T^d}{^tlD}(S_R \Theta_y U+T)ldy=I_1-I_2
\end{split}
\end{equation}
With
\begin{equation}
\begin{split}
I_1=\int_{ \T^d\times \T^d}(l-\nabla_x \chi_l(x,y))(I_d-\nabla 
\chi^{U}_.(Rx+y))(l-\nabla \chi^{D(U),T}_l(x))m_{U(R.+y)+T(.)}(dx)\,dy
\end{split}
\end{equation}
and
\begin{equation}\label{uhgdvvczzz1}
I_2=\int_{ \T^d\times \T^d}(l-\nabla_x \chi_l(x,y))\nabla 
\chi^{U}_.(Rx+y)(\nabla \chi^{D(U),T}_l(x)-l)m_{U(Rx+y)+T(.)}(dx)\,dy
\end{equation}
It is easy to see that $\chi^{D(U),T}_.$ is a minimizer in the 
variational formula \eref{eqforduvarfjhoiuh71jGer} associated to  $D(U,T)$, 
which is the effective diffusivity corresponding to two-scale 
homogenization on $U,T$ with complete separation between the scales, that is to 
say:
\begin{equation}\label{sjddskdjhv675A}
D(U,T)=\int_{x\in \T^d}{^t(}I_d-\nabla 
\chi^{D(U),T}_.(x))D(U)(I_d-\nabla \chi^{D(U),T}_.(x)) m_T(dx)
\end{equation}
A simple use of Cauchy-Schwarz inequality gives an upper bound on 
$I_1$,
\begin{equation}\label{ksjkausviuzs76521}
\begin{split}
I_1\leq& \Big(\int_{(x,y)\in (\T^d)^2}|l-\nabla_x \chi_l(x,y)|^2 
m^{U(Rx+y)+T(x)} (dx)\,dy\Big)^\frac{1}{2} \\ \times&\Big(\int_{(x,y)\in 
(\T^d)^2}|(I_d-\nabla \chi^{U}_.(Rx+y))(l-\nabla \chi^{D(U),T}_l(x))|^2 
m^{U(Rx+y)+T(x)} (dx) ,dy\Big)^\frac{1}{2}
\end{split}
\end{equation}
Integrating first in $y$ in the second term and, using the formulas 
linking effective diffusivities and solutions of the cell problem, we 
obtain
\begin{equation}\label{hgfutztvvztv1}
I_1\leq \Big(\int_{y\in \T^d}{^tlD}(S_R \Theta_y 
U+T)l\,dy\Big)^\frac{1}{2} \times \Big({^tlD}(U,T)l\Big)^\frac{1}{2}e^{\frac{\| 
T\|_\alpha}{R^\alpha}}
\end{equation}
We now estimate $I_2$. The following lemma together  with  
\eref{dsgutgveq17h71} and \eref{hgfutztvvztv1} gives lemma \ref{jhibuhzbzusbb1a}.
\begin{Lemma}\label{jhhasjssz8719}
\begin{equation}\label{eqhbibsjoanck0sh7nb63f200}
\begin{split}
|I_2| \leq &\Big(\int_{y\in \T^d}{^tlD}(S_R \Theta_y 
U+T)l\,dy\Big)^\frac{1}{2} \times \Big({^tlD}(U,T)l\Big)^\frac{1}{2} \\&C_de^{C_d 
\Osc(U)} e^{4\frac{\|T\|_\alpha}{R^\alpha}}(e^{8\frac{\|T\|_\alpha}{R^\alpha}}
-1)^\frac{1}{2}
\end{split}
\end{equation}
\end{Lemma}
The proof of this lemma relies heavily on the following elliptic type 
estimate
\begin{Lemma}\label{ksahsjvkhs8761}
\begin{equation}\label{eqcontrchi9}
\|\chi_l^U\|_\infty \leq C_d \exp\big((3d+2)\Osc(U)\big)|l|
\end{equation}
\end{Lemma}
This lemma is a consequence of   theorem 5.4,  chapter 5 of \cite{St1} 
on elliptic equations with discontinuous coefficients(see also 
\cite{St2}), we give the proof of lemma \ref{ksahsjvkhs8761} for the sake of 
completeness in paragraph \ref{hcsajscshsgcvzsvc5}.\\
We will now prove  lemma \ref{jhhasjssz8719}. First we will estimate 
the distance between $\chi_l(x,y)$ and $\chi_l(x+y/R,0)$ for $y \in 
[0,1]^d$ in $H^1$ norm.\\
By the orthogonality property of the solution of the cell problem
 for $y\in [0,1]^d$
\begin{equation}\label{Matsub2scupboeqvprjhb671AA}
\begin{split}
\int_{x \in \T^d}|\nabla_x &\chi_l(x+\frac{y}{R},0)-\nabla_x 
\chi_l(x,y)|^2\frac{e^{-2(U(Rx+y)+T(x))}}{\int_{\T^d}e^{-2(U(Rz+y)+T(z))}dz}dx\,dy\\
=&\int_{x \in \T^d}|l-\nabla_x 
\chi_l(x+\frac{y}{R},0)|^2\frac{e^{-2(U(Rx+y)+T(x))}}{\int_{\T^d}e^{-2(U(Rz+y)+T(z))}dz}dx\,dy-{^tlD}(S_R\Theta_y 
U+T)l\\
\leq & {^tlD}(S_R U+T)l 
e^{4\frac{\|T\|_\alpha}{R^\alpha}}-{^tlD}(S_R\Theta_y U+T)l
\end{split}
\end{equation}
Thus by  lemma \ref{jhibuhzbzusbb1}, for $y\in [0,1]^d$,
\begin{equation}\label{Matsub2scupboeqvprjhb671}
\begin{split}
\int_{x \in \T^d}|\nabla_x \chi_l(x+\frac{y}{R},0)-\nabla_x 
\chi_l(x,y)|^2& m_{U(R.x+y)+T(.)}(dx)\,dy
\\&\leq  {^tlD}(S_R\Theta_y U+T)l 
(e^{8\frac{\|T\|_\alpha}{R^\alpha}}-1)
\end{split}
\end{equation}
Let us introduce
\begin{equation}\label{eqhbibsjoanck0sh7nb63f1}
\begin{split}
I_3=&\int_{ (x,y)\in \T^d\times [0,1]^d}(l-\nabla_x 
\chi_l(x+\frac{y}{R},0))\nabla \chi^{U}_.(Rx+y)(\nabla 
\chi^{D(U),T}_l(x)-l)\\&\frac{e^{-2(U(Rx+y)+T(x+\frac{y}{R}))}}{\int_{\T^d}e^{-2U(z)}dz\int_{\T^d}e^{-2T(z)}dz}dx\,dy
\end{split}
\end{equation}
one has
\begin{equation}\label{eqhbibsjoanck0sh7nb63f2}
\begin{split}
|I_2-I_3| \leq &\Big(\int_{y\in \T^d}{^tlD}(S_R \Theta_y 
U+T)l\,dy\Big)^\frac{1}{2} \times \Big({^tlD}(U,T)l\Big)^\frac{1}{2}\\&6 
e^{4\frac{\|T\|_\alpha}{R^\alpha}+\Osc(U)}(e^{8\frac{\|T\|_\alpha}{R^\alpha}}
-1)^\frac{1}{2}
\end{split}
\end{equation}
This can be seen using  \eref{Matsub2scupboeqvprjhb671}, 
\eref{uhgdvvczzz1} and by Cauchy Schwarz inequality (the computation is similar to 
the one in \eref{ksjkausviuzs76521}); Voigt-Reiss' inequality ($D(U)\geq 
e^{-2\Osc(U)}$) and noticing that  for $y\in 
[0,1]^d$,$|T(x+y/R)-T(x)|\leq \|T\|_\alpha/R_\alpha$.\\
We now want to estimate $I_3$. Noting that
$$\nabla_y\Big(e^{-2(U(Rx+y)+T(x+\frac{y}{R}))}(l-\nabla_x 
\chi_l(x+\frac{y}{R},0))\Big)=0$$ $$\nabla \chi^{U}_.(Rx+y)=\nabla_y 
\chi^{U}_.(Rx+y)$$ and integrating by parts in $y$, one obtains
 \[\begin{split}
I_3=&\sum_{i=1}^d \int_{ x\in \T^d,y^i \in 
\partial^i([0,1]^d)}\Big(e^{-2T(x+\frac{y^i+e_i}{R})}(l-\nabla_x \chi_l(x+\frac{y^i+e_i}{R},0)) 
\\&-e^{-2T(x+\frac{y^i}{R})}(l-\nabla_x 
\chi_l(x+\frac{y^i}{R},0))\Big).e_i \chi^{U}_.(Rx+y^i)\\&(\nabla 
\chi^{D(U),T}_l(x)-l)\frac{e^{-2U(Rx+y^i)}}{\int_{\T^d}e^{-2U(z)}dz\int_{\T^d}e^{-2T(z)}dz}dx\,dy^i
\end{split}
\]
Where we have used the notation $$\partial^i([0,1]^d)=\{x\in 
[0,1]^d\,:\,x_i=0\}$$
Let us introduce 
\begin{equation}\label{hgccccitcc21}
\begin{split}
I_4=&\sum_{i=1}^d \int_{ x\in \T^d,y^i \in 
\partial^i([0,1]^d)}\Big(-\nabla_x \chi_l(x+\frac{y^i+e_i}{R},0) +\nabla_x 
\chi_l(x+\frac{y^i}{R},0)\Big).e_i \\&\chi^{U}_.(Rx+y^i)(\nabla 
\chi^{D(U),T}_l(x)-l)e^{-2T(x+\frac{y^i}{R})}\frac{e^{-2U(Rx+y^i)}}{\int_{\T^d}e^{-2U(z)}dz\int_{\T^d}e^{-2T(z)}dz}dx\,dy^i
\end{split}
\end{equation}
It is easy to obtain 
\begin{equation}\label{eqcontrchi6}
\begin{split}
|I_4-I_3|\leq d e^{3\Osc(U)} e^{\frac{2 \| 
T\|_\alpha}{R^\alpha}}(e^{\frac{2 \| T\|_\alpha}{R^\alpha}}-1)\Big({^tl}D(S_R 
U+T)l\Big)^\frac{1}{2}
 \|\chi^U_.\|_{\infty} \Big({^tl}D(U,T)l\Big)^\frac{1}{2}
\end{split}
\end{equation}
We will now put into evidence the fact that although $\chi_l(x,0)$ is 
not periodic on $R^{-1}\T^d$ the distance  (with respect to the natural 
$H^1$ norm) between the solution of the cell problem $\chi_l(x,0)$ and 
its translation $\chi_l(x+e_k/R,0)$ along the axis of the torus 
$R^{-1}\T^d$,  is small.  This is due to the presence of a fast period $R^{-1} 
\T^d$ in the decomposition $V=S_R U+T$.\\
Using the standard property of the solution of the cell problem one 
obtains
\begin{equation}
\begin{split}
\int_{\T^d}&|\nabla \chi_l (x+\frac{e_k}{R},0)-\nabla \chi_l (x,0)|^2 
m_{S_R U+T}(dx)\\
&=\int_{\T^d} |l-\nabla \chi_l (x,0)+\nabla \chi_l 
(x+\frac{e_k}{R},0)-\nabla \chi_l (x,0)|^2 m_{S_R U+T}(dx)
- {^tl}D(S_RU+T)l
\\&\leq e^{4\frac{\|T\|_\alpha}{R^\alpha}} \int_{\T^d} |l-\nabla \chi_l 
(x+\frac{e_k}{R},0))|^2 m_{\Theta{\frac{e_k}{R}}(S_R U+T)}(dx)
- {^tl}D(S_RU+T)l
\end{split}
\end{equation}
which leads to 
\begin{equation}\label{eqcontrchi7}
\int_{\T^d}|\nabla \chi_l (x+\frac{e_k}{R},0)-\nabla \chi_l (x,0)|^2 
m^{S_R U+T}(dx) \leq 
{^tl}D(S_RU+T)l(e^{4\frac{\|T\|_\alpha}{R^\alpha}}-1)
\end{equation}
Combining this inequality with the definition \eref{hgccccitcc21} of 
$I_4$, 
Cauchy-Schwarz inequality proves that
\begin{equation}
\begin{split}
|I_4|\leq &\sum_{i=1}^d \int_{y^i \in \partial^i([0,1]^d)} \Big(\int_{ 
x\in \T^d} \big((-\nabla_x \chi_l(x+\frac{y^i+e_i}{R},0) \\&+\nabla_x 
\chi_l(x+\frac{y^i}{R},0)) .e_i\big)^2 
\frac{e^{-2U(Rx+y^i)-2T(x+\frac{y^i}{R})}}{\int_{\T^d}e^{-2U(z)}dz\int_{\T^d}e^{-2T(z)}dz}dx 
\Big)^\frac{1}{2}\\
\Big(&\int_{ x\in \T^d} \big(\chi^{U}_.(Rx+y^i)(\nabla 
\chi^{D(U),T}_l(x)-l)\big)^2 
\frac{e^{-2U(Rx+y^i)-2T(x+\frac{y^i}{R})}}{\int_{\T^d}e^{-2U(z)}dz\int_{\T^d}e^{-2T(z)}dz}dx  \Big)^\frac{1}{2}\,dy^i
\end{split}
\end{equation}
Combining this with \eref{eqcontrchi7} one obtains
\begin{equation}\label{eqcontrchi8}
\begin{split}
|I_4|\leq (e^{\frac{8 \| 
T\|_\alpha}{R^\alpha}}-1)^\frac{1}{2}\Big({^tlD}(S_R U+T)l\Big)^\frac{1}{2}
d 
\|\chi^U_.\|_{\infty}e^{3\Osc(U)}\Big({^tl}D(U,T)l\Big)^\frac{1}{2}e^{\frac{2 \| T\|_\alpha}{R^\alpha}}
\end{split}
\end{equation}
Using  lemma \ref{eqcontrchi9} to estimate $\|\chi^U_.\|_{\infty}$ in  
\eref{eqcontrchi8} and 
combining \eref{eqcontrchi6} and \eref{eqhbibsjoanck0sh7nb63f2} one 
obtains \eref{eqhbibsjoanck0sh7nb63f200} and  lemma \ref{jhhasjssz8719}, 
which proves the upper bound of theorem \ref{NewthcontmscM}.

\paragraph{}\label{hcsajscshsgcvzsvc5}
The purpose of this paragraph is to prove the estimate 
\eref{eqcontrchi9}. First we will remind a theorem concerning elliptic equation with 
discontinuous coefficient from G. Stampacchia. Its proof in a more 
general form can be found in \cite{St2}, chapter 5, theorem 5.4 (see also 
\cite{St1}).\\
Let us consider the operator (in the weak sense)$L=\nabla(A\nabla)$
defined on some open set $\Omega\subset \R^d$ (for $d\geq 3$) with 
smooth boundary $\partial \Omega$. $A$ is a $d \times d$ matrix with 
bounded  coefficients in $L^\infty(\Omega)$ such that for all $\xi\in\R^d$, 
$\lambda |\xi|^2 \leq {^t\xi}A\xi$ 
and for all $i,j$; $|A_{ij}|\leq M$, for some positive constant 
$0<\lambda, M<\infty$.\\
Let  $p>d\geq 3$. 
For $1 \leq i \leq d$ let $f_i \in L^p(\Omega)$\\
if $\chi \in H^1_{loc}(\Omega)$ is a local (weak) solution of the 
equation
\begin{equation}\label{AnToDivOpAnafipar}
\nabla \big(A\nabla \chi\big)=-\sum\limits_{i=1}^{d}\partial_i f_i 
\end{equation}
 then $\chi$ is in $L^\infty(\Omega)$ and if $x_0 \in \Omega$ and $R>0$
\begin{Theorem}\label{Stampa1}
The solution of \eref{AnToDivOpAnafipar} verify the following 
inequality (in the essential supremum sense with $\Omega(x_0,R)=\Omega \cap 
B(x_0,R)$) 
\begin{equation}
\begin{split}
\max_{\Omega(x_0,\frac{R}{2})} |\chi| \leq K \Big[ \big\{\frac{1}{R^d} 
\int_{\Omega(x_0,R)} \|\chi\|^2 \big\}^{\frac{1}{2}} 
 +\sum\limits_{i=1}^{d} \|f_i\|_{L^p(\Omega(x_0,R))} 
\frac{R^{1-\frac{d}{p}}}{\lambda}  \Big]
\end{split}
\end{equation}
with $K=C_d (\frac{M}{\lambda})^\frac{3d}{2}$
\end{Theorem}
The explicit dependence of the constants in $M$ and $\lambda$ have been 
obtained by following the proof of G. Stampacchia \cite{St2}.
We will now prove \eref{eqcontrchi9} for $d\geq 3$ (For $d=1$, this 
estimate is trivial, for $d=2$, it is sufficient consider $U(x_1,x_2)$ as 
a function on $T^3_1$ to obtain the result).
$\chi_l$ satisfies 
\[\nabla \big(\exp(-2U)\nabla \chi_l\big)=l.\nabla \exp(-2U)\]
then by theorem \ref{Stampa1}  for $x_0 \in [0,1]^d$
\[\max_{B(x_0,\frac{1}{2})}|\chi_l| \leq C_d \exp(3\Osc(U)d) \big[ 
\big( \int_{B(x_0,1)} |\chi_l|^2 \big)^{\frac{1}{2}} 
+|l| \exp(2 \Osc(U)) \big]\]
Now by periodicity
\[\int_{B(x_0,1)} |\chi_l|^2 dx \leq \int_{\T^d} |\chi_l|^2 dx\]
and by  Poincar\'{e} inequality (we  assume $\int_{\T^d}\chi_l(x)dx=0$)
\[\int_{\T^d} |\chi_l|^2 dx \leq C_d \int_{\T^d} |\nabla \chi_l|^2 dx 
\]
thus 
\[\int_{B(x_0,1)} |\chi_l|^2 dx \leq C_d \exp(2\Osc(U))  \int_{\T^d} 
|\nabla \chi_l|^2 m_U(dx)\]
And since
\[\int_{\T^d} |l-\nabla \chi_l|^2 m_U(dx)= l^2-\int_{\T^d} |\nabla 
\chi_l|^2 m_U(dx)\]
one has
\[\int_{\T^d} |\nabla \chi_l|^2 m_U(dx) \leq l^2\]
and the bound on  $\|\chi_l\|_\infty$ is proven.

\subsubsection{Quantitative multi-scale-homogenization: Lower bound in 
the  theorem \ref{NewthcontmscM}}
\paragraph{}\label{hcsajscsddddaavzsvc1}
As for the upper bound it is sufficient to prove theorem  
\ref{NewthcontmscM}
assuming that $U$ and $T$ are smooth and we will use the notation 
introduced in the paragraph \ref{hcsajscshsgcvzsvc1}.
We will prove below that
\begin{Lemma}\label{djksdhhhhhh86681}
If $R\geq \|T\|_\alpha$ then for $\xi \in \S^{d-1}$
\begin{equation}
\int_{\T^d}{^t\xi}D(S_R U+\Theta_y T)^{-1}\xi dy\leq 
(1+C_{d,\Osc(U),\|U\|_\alpha,\alpha,\|T\|_\alpha} R^{-\alpha/2}){^t\xi}D(U,T)^{-1}\xi
\end{equation}
\end{Lemma}
This lemma with lemma \ref{jhibuhzbzusbb1} gives the lower bound in  
theorem \ref{NewthcontmscM}. We now prove  lemma \ref{djksdhhhhhh86681}.
Let us introduce
\begin{equation}\label{eqexplfornak81}
P(x,y) =I_d-\frac{\exp(-2(S_R \Theta_y U+T))}{\int_{\T^d}\exp(-2(S_R 
\Theta_y U+T)(x))dx} (I_d-\nabla \chi(x,y)_.)D(S_R \Theta_y U+T)^{-1}
\end{equation}
\begin{equation}\label{eqexplforna5k812}
P^U(x) =I_d-\frac{\exp(-2 U(x))}{\int_{\T^d}\exp(-2U(x))dx} (I_d-\nabla 
\chi^U(x)_.)D(U)^{-1}
\end{equation}
 and
\begin{equation}
P^{D(U),T}(x)=I_d-\frac{e^{-2T(x)}}{\int_{\T^d}e^{-2T(x)}dx} 
D(U)(I_d-\nabla \chi^{D(U),T}(x)_.)D(U,T)^{-1}
\end{equation}
We remind that $P(x,y)$ minimize the well known variational formula 
associated to $D(S_R \Theta_y U+T)^{-1}$, that is why it will play for the 
lower bound in the  theorem \ref{NewthcontmscM} the role played by the 
gradient of the solution of the solution of the cell problem  $\nabla 
\chi(x,y)_.$ for the upper bound. More precisely,  for $\xi \in 
\S^{d-1}$ one obtains as in the proof of the upper bound (by decomposing $\xi$ 
here)
\begin{equation}\label{eqkjhahhsb7h7283}
\begin{split}
\int_{y\in \T^d}&{^t\xi}D(S_R \Theta_y U+T)^{-1}\xi \,dy 
\\&=\int_{(x,y)\in (\T^d)^2}\Big(\int_{\T^d} e^{-2(S_R \Theta_y 
U+T)(z)}dz\Big)e^{2(S_R \Theta_y U+T)(x)}{^t\xi}{^t(}I_d-P(x,y))\xi\,dx\,dy\\
&\leq e^{\frac{2\| T\|_\alpha}{R^\alpha}} (I_1+I_2)
\end{split}
\end{equation}
with
\begin{equation}
\begin{split}
I_1=&\int_{\T^d}e^{-2U(z)}\,dz \int_{\T^d}e^{-2T(z)}\,dz
\int_{(x,y)\in (\T^d)^2} e^{2(S_R \Theta_y 
U+T)(x)}{^t\xi}{^t(}I_d-P(x,y))\\&(I_d-P^U(Rx+y))(I_d-P^{D(U),T}(x))\xi\,dx\,dy
\end{split}
\end{equation}
and
\begin{equation}
\begin{split}
I_2=&\int_{\T^d}e^{-2U(z)}\,dz \int_{\T^d}e^{-2T(z)}\,dz
\int_{(x,y)\in (\T^d)^2} e^{2(S_R \Theta_y 
U+T)(x)}{^t\xi}{^t(}I_d-P(x,y))\\&P^U(Rx+y)(I_d-P^{D(U),T}(x))\xi\,dx\,dy
\end{split}
\end{equation}
As for the upper bound, using Cauchy Schwarz inequality for the 
integration in $x$ and $y$, and using
\begin{equation}
\begin{split}
{^t\xi}D(U,T)^{-1}\xi=&\int_{(x,y)\in (\T^d)^2} e^{2(S_R \Theta_y 
U+T)(x)}\\&\big((I_d-P^U(Rx+y))(I_d-P^{D(U),T}(x))\xi\big)^2 dx\,dy
\end{split}
\end{equation}
one obtains that
\begin{equation}\label{eqkjhahhsb7h7282}
\begin{split}
|I_1|\leq  e^{\frac{\| T\|_\alpha}{R^\alpha}}\Big(\int_{y\in 
\T^d}{^t\xi}D(S_R \Theta_y U+T)^{-1}\xi \,dy \Big)^\frac{1}{2}
\Big({^t\xi}D(U,T)^{-1}\xi \Big)^\frac{1}{2}
\end{split}
\end{equation}
Thus $I_2$ will be an error term and it will be proven below that
\begin{Lemma}\label{jhhhhhhsiudz16}
\begin{equation}\label{eqkjhahhsb7h7281}
\begin{split}
|I_2|\leq  &\Big(\int_{y\in \T^d}{^t\xi}D(S_R \Theta_y U+T)^{-1}\xi 
\,dy \Big)^\frac{1}{2}
\Big({^t\xi}D(U,T)^{-1}\xi 
\Big)^\frac{1}{2}\\&C_{d,\Osc(U),\|U\|_\alpha,\alpha} e^{4\frac{\| T\|_\alpha}{R^\alpha}}(e^{8\frac{\| 
T\|_\alpha}{R^\alpha}}-1)^\frac{1}{2}
\end{split}
\end{equation}
\end{Lemma}
Let us observe that combining the estimate \eref{eqkjhahhsb7h7281} of 
lemma \ref{jhhhhhhsiudz16} with \eref{eqkjhahhsb7h7282} and  
\eref{eqkjhahhsb7h7283} proves lemma \ref{djksdhhhhhh86681}.\\
We will now prove the  lemma \ref{jhhhhhhsiudz16}.
As it has been done in the proof of the upper bound it is easy to show 
that, with
\begin{equation}\label{eqkjhahhsb7h7280033}
\begin{split}
I_3=&\int_{\T^d}e^{-2U(z)}\,dz \int_{\T^d}e^{-2T(z)}\,dz
\int_{(x,y)\in \T^d\times [0,1]^d} 
e^{2(U(Rx+y)+T(x+\frac{y}{R}))}\\&{^t\xi}{^t(}I_d-P(x+\frac{y}{R},0))P^U(Rx+y)(I_d-P^{D(U),T}(x))\xi\,dx\,dy
\end{split}
\end{equation}
one has
\begin{equation}\label{eqkjhahhsb7h7280032}
\begin{split}
|I_3-I_2|\leq 6 e^{\Osc(U)}e^{\frac{4\| 
T\|_\alpha}{R^\alpha}}(e^{8\frac{\| 
T\|_\alpha}{R^\alpha}}-1)^\frac{1}{2}\Big({^t\xi}D(S_RU+T)\xi\Big)^\frac{1}{2}\Big({^t\xi}D(U,T)\xi\Big)^\frac{1}{2}
\end{split}
\end{equation}
It will be proven in \ref{hcsajscsddddaavzsvc2} that
\begin{Lemma}\label{Matsubskewsymassokhulkjh1}
There exists $d\times d\times d$ tensors $H_{ijm}^U$ such that 
$H_{ijm}^U=-H_{jim}^U\in C^\infty(\T^d)$,
\begin{equation}\label{eqkjhahhsb7h7280034}
P_{im}^U=\sum_{j=1}^d \partial_j H_{ijm}^U\quad \text{and}\quad 
\|H_{ijm}^U\|_\infty \leq C_{d,\Osc(U),\|U\|_\alpha,\alpha}
\end{equation}
\end{Lemma}
Combining \eref{eqkjhahhsb7h7280034} with the explicit formula 
\eref{eqexplfornak81} for $P$ one obtains
\begin{equation}
\begin{split}
I_3=&\frac{\int_{\T^d}e^{-2U(z)}\,dz 
\int_{\T^d}e^{-2T(z)}}{\int_{\T^d}\exp(-2(S_R U+T)(z))dz} \int_{(x,y)\in \T^d\times 
[0,1]^d}\sum_{i,j,k=1}^d({^t\xi}{^t(}I_d-P^{D(U),T}(x)))_i\\&
\partial_k H^{U}_{i,k,j}(Rx+y)\Big((I_d-\nabla 
\chi_.(x+\frac{y}{R})\big)D(S_R  U+T)^{-1}\xi\Big)_j \,dx\,dy
\end{split}
\end{equation}
Thus, using the same notation as in the equation \eref{hgccccitcc21} 
and integrating by parts in $y$, one obtains
\begin{equation}
\begin{split}
I_3&=\frac{\int_{\T^d}e^{-2U(z)}\,dz 
\int_{\T^d}e^{-2T(z)}}{\int_{\T^d}\exp(-2(S_R U+T)(z))dz} \sum_{i,j,k=1}^d \int_{(x,y^k)\in \T^d\times 
\partial^k([0,1]^d)}({^t\xi}{^t(}I_d-P^{D(U),T}(x)))_i\\&
 H^{U}_{i,k,j}(Rx+y^k)\Big((\nabla \chi_.(x+\frac{y^k}{R})-\nabla 
\chi_.(x+\frac{y^k+e_k}{R})\big)D(S_R  U+T)^{-1}\xi\Big)_j \,dx\,dy^k
\end{split}
\end{equation}
Which, using Cauchy Schwarz inequality, leads to
\begin{equation}
\begin{split}
|I_3|\leq& C_d e^{2\Osc(U)}\sup_{ijk}\|H^{U}_{i,k,j}\|_\infty  
\sum_{k=1}^d \int_{y^k \in  \partial^k([0,1]^d)}   \Big(\int_{x\in 
\T^d}\big((I_d-P^{D(U),T}(x))\xi\big)^2 \\& e^{2(U(Rx+y^k)+T(x+\frac{y^k}{R}))} 
dx\Big)^\frac{1}{2}
\Big(\int_{x\in \T^d}\Big(\big((\nabla \chi_.(x+\frac{y^k}{R})-\nabla 
\chi_.(x+\frac{y^k+e_k}{R})\big)\\&D(S_R  U+T)^{-1}\xi\Big)^2 
e^{-2(U(Rx+y^k)+T(x+\frac{y^k}{R}))} dx\Big)^\frac{1}{2} \,dy^k
\end{split}
\end{equation}
Using  bounds \eref{eqkjhahhsb7h7280034}  and the equation 
\eref{eqcontrchi7} to control the natural $H^1$ distance between the  solution of 
the cell problem $\chi_.(x+y^k/R)$ and its translation by $e_k/R$ one 
obtains
\begin{equation}\label{weuewubobaeq191}
\begin{split}
|I_3|\leq& \Big(\xi D(S_R  U+T)^{-1}\xi\Big)^\frac{1}{2}\Big({^t\xi 
D}(U,T)^{-1}\xi\Big)^\frac{1}{2}\\&C_{d,\Osc(U),\| U\|_\alpha,\alpha} 
e^{4\frac{\| T\|_\alpha}{R^\alpha}} (e^{4\frac{\| 
T\|_\alpha}{R^\alpha}}-1)^\frac{1}{2}
\end{split}
\end{equation}
Combining \eref{weuewubobaeq191} and \eref{eqkjhahhsb7h7280032}  one 
obtains \eref{eqkjhahhsb7h7281}, which proves lemma 
\ref{eqkjhahhsb7h7281}.

\paragraph{}\label{hcsajscsddddaavzsvc2}
In this paragraph, we will prove lemma \ref{Matsubskewsymassokhulkjh1}.
Since for each $m\in \{1,\ldots,d\}$, $P_{.,m}^U$ are divergence free 
vectors with mean $0$ with respect to Lebesgue measure, by the 
proposition 4.1 of \cite{JiKo99} there exists a skew-symmetric $\T^d$-periodic 
smooth matrices $H_{ij1}^U,\ldots,H_{ijd}^U$ ($H_{ijm}^U=-H_{jim}^U$) 
such that for all $m$
\begin{equation}
P_{im}^U=\sum_{j=1}^d \partial_j H_{ijm}^U
\end{equation}
Moreover writing
\begin{equation}
P_{.m}^U=\sum_{k\not= 0}p^k_{.m} e^{2i\pi (k.x)}
\end{equation}
the Fourier series expansion of $P^U$, one has (see the proposition 4.1 
of \cite{JiKo99})
\begin{equation}
H_{njm}^U=\frac{1}{2i\pi} \sum_{k\not= 0} 
\frac{p^k_{nm}k_j-p^k_{jm}k_n}{k^2} e^{2i\pi (k.x)}
\end{equation}
Let us observe that
\begin{equation}\label{eqexplforna5k812shn1}
H_{njm}^U=B_{nm}^j-B_{jm}^n
\end{equation}
where $B_{nm}^j$ are the smooth $\T^d$-periodic solutions of  $\Delta 
B_{nm}^j=\partial_j P_{nm}^U$.
By theorem \ref{Stampa1} (theorem 5.4,  chapter 5 of \cite{St1}),  if 
$B_{nm}$ is chosen so that $\int_{\T^d}B_{nm}(x)dx=0$ then 
$\|B_{nm}^j\|_\infty \leq C_d \|P_{nm}^U\|_\infty$. Now using theorem 1.1 of 
\cite{YaVo00}  it is easy to obtain that
$\|\nabla \chi^U_l\|_\infty\leq C_{d,\Osc(U),\|U\|_\alpha,\alpha}$, 
combining this with \eref{eqexplforna5k812}, one obtains
\begin{equation}
\|B_{nm}^j\|_\infty \leq C_{d,\Osc(U),\|U\|_\alpha,\alpha}
\end{equation}
Which leads to \eref{eqkjhahhsb7h7280034} by the equation 
\eref{eqexplforna5k812shn1}.

\subsubsection{Explicit formulas of ef{f}ective dif{f}usivities from  
level-3 large deviations. Proof of theorem \ref{assjhsao78982hh1}}
The equation \eref{dssdddjsdj5365} follows from the Voigt-Reiss 
inequality: for $U\in L^\infty(\T^d)$
\begin{equation}
D(U)\geq I_d \big(\int_{\T^d}e^{2 U(x)}\int_{\T^d}e^{-2 U(x)}\big)^{-1}
\end{equation}
and the fact that, if $U\in C^\alpha(\T^d)$, then by the Varadhan's 
lemma and level-3 large deviation associated to the shift $s_\rho$,
\begin{equation}
\lim_{n\rightarrow \infty}\frac{1}{n} \ln 
\Big(\int_{\T^d}e^{\sum_{k=0}^{n-1} U(s_\rho^k x)}dx\Big)=P_\rho(U) 
\end{equation}
We refer to \cite{Ow00a} for a more detailed proof of this statement.\\
In higher dimensions, when the medium is self-similar one can use the 
criterion \eref{eqjhl9889b87k2} associated with the equation 
\eref{dssdddjsdj5365} to characterize ratios for which $D(V_0^n)$ does not 
converge to $0$ with an exponential rate.
The equation \eref{assssskavuuszz}, i.e. the extension of the result 
\eref{eqjhl9889b87k1} to dimension $2$ is done by observing that
\begin{Proposition}\label{Thjbztvjhvzu7h71}
For $d=2$ one has
\begin{equation}\label{eqahsdihdbbodbbbz71}
\begin{split}
\lambda_{\max}\big(D(U)\big)\lambda_{\min}\big(D(-U)\big)&=\lambda_{\min}\big(D(U)\big)\lambda_{\max}\big(D(-U)\big)\\
&=\frac{1}{\int_{T^d_1} \exp(2U)dx \int_{T^d_1}\exp(-2U) dx}
\end{split}
\end{equation}
\end{Proposition}
\noindent From which one deduces that if $D(U)=D(-U)$ then 
\begin{equation}
\lambda_{\min}(D(U))\lambda_{\max}(D(U))=\big(\int_{T^d_1} \exp(2U)dx 
\int_{T^d_1}\exp(-2U) dx\big)^{-1}
\end{equation}
Let us observe that the assumption $D(U)=D(-U)$ is satisfied if, for 
instance $-U_n(x)=U_n(-x)$ or $-U_n(x)=U_n(Ax)$ where $A$ is an isometry 
of $\R^d$. And the existence of a reflection $B$ such that $U(Bx)=U(x)$ 
ensures that $\lambda_{\min}(D(U))=\lambda_{\max}(D(U))$. Thus
these symmetry hypotheses combined with \eref{assssskavuuszz} ensure 
the validity of \eref{assssskavuuszz212}.\\
It would be interesting to extend the equation \eref{assssskavuuszz} of 
theorem \ref{assjhsao78982hh1} to more general cases and higher 
dimensions. Indeed the proposition \ref{Thjbztvjhvzu7h71} is deduced from the 
following proposition \ref{Matsubcodqlamjhga2} that put into evidence a 
strong geometrical link between cohomology and homogenization.\\
Write $\F_{sol}=\Big\{p \in (C^\infty(T^d_1))^d | 
\operatorname{div}(p)=0 \;  \text{and} \;  \int_{T^d_1}p dx=0 \Big\}$
and $Q(U)$ the positive, definite, symmetric matrix associated to the 
following variational problem. For $l\in \S^d$
\begin{equation}\label{gvvvvzziubveq1}
^tlQ(U)l=\inf_{p\in \F_{sol}}\frac{\int_{T^d_1}|l-p|^2 
\exp(2U)dx}{\int_{T^d_1} \exp(2U)dx}
\end{equation}
Write in the increasing order $\lambda(D(U))_i$ and decreasing order 
$\lambda(Q(U))_i$ the eigenvalues of $D(U)$ and $Q(U)$.
\begin{Proposition}\label{Matsubcodqlamjhga2}
For all $i\in \{1,\ldots,d\}$
\begin{equation}\label{Matrixsubeqdqugv3}
\lambda(D(U))_i \lambda(Q(U))_i= \frac{1}{\int_{T^d_1} \exp(2U)dx 
\int_{T^d_1}\exp(-2U) dx}
\end{equation}
\end{Proposition}
Now we will introduce a geometric interpretation of homogenization that 
will allow us to prove proposition \ref{Thjbztvjhvzu7h71} and equation 
\eref{assssskavuuszz} of theorem \ref{assjhsao78982hh1}.
Let $U \in C^{\infty}(\T^d)$. It is easy to obtain the following 
orthogonal decomposition
\begin{equation}\label{matsuborthdecompH}
H=(L^2(m_U))^d=H_{pot}\oplus H_{sol}
\end{equation}
Where $H_{pot}$, $H_{sol}$ are the closure (with respect to the 
intrinsic norm $\|.\|_H$) of the sets of $C_{pot},C_{sol}$  the sets of 
smooth, $\T^d$-periodic, potential and solenoidal vector fields, i.e. with 
$\C=(C^\infty(\T^d))^d$
\begin{equation}
\C_{pot}=\Big\{\xi \in \C \,|\, \exists f \in C^{\infty}(\T^d) \;  
\text{with} \;  \xi=\nabla f \Big\}
\end{equation}
\begin{equation}
\C_{sol}=\Big\{\xi \in \C | \exists p \in \C \;  \text{with} \;  
\operatorname{div}(p)=0 \;  \text{and} \;  \xi=p \exp(2U)\int_{\T^d} 
e^{-2U(x)}dx\Big\}
\end{equation}
Thus H is a real Hilbert space equipped with the scalar product\\ 
$(\xi,\nu)_H=\int_{\T^d}\xi(x).\nu(x)m_U(dx)$
and by the variational formulation \eref{eqforduvarfjhoiuh71}, for $l 
\in \R^d$, $^tlD(U)l$ is the norm in $H$ of the orthogonal projection of 
$l$ on $H_{sol}$ and $l=\nabla \chi_l + \exp(2U) p_l$ is the orthogonal 
decomposition of $l$.
\begin{equation}
\sqrt{^tlD(U)l}=\dist(l,H_{pot})
\end{equation}
Now by duality for all $\xi \in H$
\begin{equation}
\dist(\xi,H_{pot})= \sup_{\delta \in \C_{sol}} 
\frac{(\delta,\xi)_H}{\|\delta\|_H}
\end{equation}
From which we deduce  the following variational formula for the 
effective diffusivity by choosing $\xi=l \in \R^d$
\begin{equation}\label{Matrixsubeqdqugv1}
^tlD(U)l=
\sup_{p\in \C \; \operatorname{div}(p)=0} \frac{\big(\int_{\T^d}l.p 
dx\big)^2}{\int_{\T^d}p^2 \exp(2U)dx \int_{\T^d}\exp(-2U) dx}
\end{equation}
Note  that the equation \eref{Matrixsubeqdqugv1} gives back 
Voigt-Reiss's inequality by choosing $p=l$.\\
Let $Q(U)$ be the positive, definite, symmetric matrix given by the  
variational formula \eref{gvvvvzziubveq1}.
Then the following proposition is a direct consequence of  the equation 
\eref{Matrixsubeqdqugv1}.
\begin{Proposition}\label{Matsubcodqlamjhga1}
For all $l\in \S^{d-1}$
\begin{equation}\label{Matrixsubeqdqugv2}
^tlD(U)l= \frac{1}{\int_{\T^d} \exp(2U)dx \int_{\T^d}\exp(-2U) 
dx}\sup_{\xi \in \S^{d-1}}\frac{(l.\xi)^2}{^t\xi Q(U)\xi}
\end{equation}
\end{Proposition}
Choosing an orthonormal basis diagonalizing $Q(U)$, it is an easy 
exercise to use this proposition in order to establish a one to one 
correspondence between the eigenvalues  of $Q(U)$ and $D(U)$ to obtain the 
proposition \ref{Matsubcodqlamjhga2}.

\paragraph{Dimension two}
In dimension two, the Poincar\'{e} duality establishes a simple 
correspondence between $Q(U)$ and $D(-U)$.
\begin{Proposition}\label{MatSubpropdim2kzf1}
For $d=2$, one has
\begin{equation}
Q(U)={^tP}D(-U)P
\end{equation}
where $P$ stands for the rotation matrix
\begin{equation}
P=\begin{pmatrix}
0 & -1 \\
1 & 0
\end{pmatrix}
\end{equation}
\end{Proposition}
Indeed by the Poincar\'{e} duality one has $\F_{sol}=\{P\nabla f\,:\, 
f\in C^\infty(\T^d)\}$
and the proposition \ref{MatSubpropdim2kzf1} follows from the 
definition of $Q(U)$.
The proposition \ref{Thjbztvjhvzu7h71} is then a direct consequence of 
the proposition \ref{MatSubpropdim2kzf1} and one deduces from the 
equation \eref{eqahsdihdbbodbbbz71} that if $D(U)=D(-U)$ then
\begin{equation}
\lambda_{\max}(D(U))\lambda_{\min}(D(U))=\Big(\int_{\T^d} \exp(2U)dx 
\int_{\T^d}\exp(-2U) dx\Big) I_d
\end{equation}
Which leads to the equation \eref{assssskavuuszz} of theorem 
\ref{assjhsao78982hh1} by theorem 3.1 of \cite{Ow00a}.

\subsection{Sub-diffusive behavior from homogenization on infinitely 
many scales}\label{kvdiudd22}
\subsubsection{Anomalous behavior of the exit times: Theorems 
\ref{IntSMAldsubanhianmethboaj1} and \ref{IntSMAldcontrol_taudx2}.}
\paragraph{}\label{paerxtihshsns7h51}
In this subsection we will prove the asymptotic anomalous behavior of 
the mean exit times $\E_x[\tau(0,r)]$ defined as weak solutions of $L_V 
f=-1$ with Dirichlet conditions on $\partial B(0,r)$. Here $U_n \in 
C^1(\T^d)$, nevertheless we will assume first that those functions are 
smooth and prove quantitative anomalous estimates on $\E_x[\tau(0,r)]$ 
depending only on the values of $D(V_0^n)$, $K_0$ and $K_\alpha$. Then, 
using standard estimates on the Green functions
associated to divergence form elliptic operators (see for instance 
\cite{St2}) it is easy to check that the exit times $\E_x[\tau(0,r)]$ are 
continuous with respect to a perturbation of $V$ in 
$L^\infty(B(0,r))$-norm. Using the density of smooth functions on $\overline{B(0,r)}$ in 
the set of bounded functions, we will then deduce that our estimates are 
valid for $U_n \in C^\alpha(\T^d)$.\\
Thus we can see the exit times as those associated to the solution of 
\eref{IntModelsubdiffstochdiffequ} and take advantage of the Ito 
formula.\\
The central lemma of the proof is lemma \ref{ksajhhdzz8717}, which will 
be proven in the paragraph \ref{paerxtihshsns7h52}.\\ 
Writing $m^{r}_U(dx)=e^{-2U(x)}\,dx(\int_{B(0,r)}e^{-2U(x)}\,dx)^{-1}$, 
we will prove in the paragraph
\ref{paerxtihshsns7h53} that for $P\in C^\infty(\overline{B(0,r})$
\begin{equation}\label{eq6g26hjsnm81}
\begin{split}
 \int_{B(0,r)} \E_x^{U+P}\big[\tau(0,r)\big] m^{r}_{U+P}(dx) &\leq e^{2 
\Osc_{B(0,r)}(P)} \int_{B(0,r)} \E_x^{U}\big[\tau(0,r)\big] 
m^{r}_{U+P}(dx)\\
&\geq e^{-2 \Osc_{B(0,r)}(P)} \int_{B(0,r)} 
\E_x^{U}\big[\tau(B(0,r))\big] m^{r}_{U+P}(dx)
\end{split}
\end{equation}
We give here the outline of the proof (see \cite{Ow00a} for $d=1$).  A 
perpetual homogenization process takes place over the infinite number 
of scales $0,\ldots,n,\ldots$ and the idea is still to distinguish, when 
one tries to estimate \eref{hgvscdklvxkj1},  the smaller scales  which 
have already been homogenized ($0,\ldots,n_{ef}$ called effective 
scales), the bigger scales  which have not had a visible influence on the 
diffusion ($n_{dri},\ldots,\infty$ called drift scales because they will 
be replaced by a constant drift in the proof) and some intermediate 
scales that manifest the particular geometric structure of their 
associated potentials in the behavior of the diffusion 
($n_{ef}+1,\ldots,n_{dri}-1=n_{ef}+n_{per}$ called perturbation scales because they will enter 
in the proof as a perturbation of the homogenization process over the 
smaller scales).\\
We will now use \eref{Aldtau_alphaeq3} and \eref{eq6g26hjsnm81} to 
prove theorem \ref{IntSMAldsubanhianmethboaj1}. For that purpose, we will 
first fix the number of scales that one can consider as 
\emph{homogenized} (we write $ef$ for \emph{effective})
\begin{equation}
n_{ef}(r)=\sup\{n \geq 0 \;:\;  e^{(n+1)(9d+15)K_0} R^2_n  \leq 
C_1/(8C_d) r^2\}<\infty
\end{equation}
where $C^1$ and $C_d$ are the constants appearing in the left term of 
\eref{Aldtau_alphaeq3}, next we fix the number of scales that
will enter in the computation as a perturbation of the homogenization 
process (we write $per$ for \emph{perturbation})
\begin{equation}
n_{per}(r)=\inf\{n \geq 0 \;:\; R_{n+1} \geq r\}-n_{ef}(r)
\end{equation}
For $r>C_{d,K_0,\rho_{\max}}$, $n_{ef}(r)$ and $n_{per}(r)$ are well 
defined. Let us choose $U=V_0^{n_{ef}(r)}$, $P=V_{n_{ef}(r)+1}^\infty$ in 
\eref{eq6g26hjsnm81}, we will bound from above, 
$\Osc_{B(0,r)}(V_{n_{ef}(r)+1}^\infty)$ by 
$$\Osc(V_{n_{ef}(r)+1}^{n_{ef}(r)+n_{per}(r)})+\|V_{n_{ef}(r)+n_{per}(r)+1}^\infty\|_\alpha r^\alpha$$ In the lower bound 
of \eref{eq6g26hjsnm81} when $x\in B(0,r/2)$ we will bound 
$\E_x^U[\tau(0,r)]$ from below by $\E_x^U[\tau(x,r/2)]$ and in the upper bound when 
$x\in B(0,r)$ we will bound it from above by
$\E_x^U[\tau(x,2r)]$. Then using  \eref{Aldtau_alphaeq3} to control 
those exit times
one obtains
\begin{equation}
\begin{split}
\int_{B(0,r)} \E_x^V \big[\tau(B(0,r))\big] m^{B(0,r)}_{V}(dx) &\leq 
C_d  e^{C_{K_\alpha,\alpha}+8 n_{per}(r)K_0} 
\frac{r^2}{\lambda_{\max}\big(D(V^{0,n_{ef}(r)})\big)}
 \\
&\geq C_{d} e^{-C_{K_\alpha,\alpha}- 8n_{per}(r)K_0}  
\frac{r^2}{\lambda_{\max}\big(D(V^{0,n_{ef}(r)})\big)}
\end{split}
\end{equation}
Theorem \ref{IntSMAldsubanhianmethboaj1} follows directly from the last 
inequalities by using the estimates \eref{eqesdvo0n61} on $D(V_0^n)$, 
\eref{Modsubboundrnrhonmin} on $R_n$ and observing that
\begin{equation}\label{aldanalpha_rinepr1}
n_{per}(r) \leq \inf\{m\geq 0\,:\, 
\frac{R_{m+n_{ef}(r)+1}}{R_{n_{ef}(r)+1}}\geq C_d e^{(n_{ef}(r)+2)(9d+15)K_0/2} \}
\end{equation}
The proof of theorem \ref{IntSMAldcontrol_taudx2} follows similar 
lines, the stability result \eref{eq6g26hjsnm81} being replaced by the 
stability condition \ref{IntSubModCoStCojh1}.

\paragraph{}\label{paerxtihshsns7h52}
It is sufficient to prove the equation \eref{Aldtau_alphaeq3} for 
$x=0$. \\
Let for $l\in \S^{d-1}$, $\chi_l$ be the $\T^d_R$-periodic solution of 
the cell problem associated to $L_U$ with $\chi_l(0)=0$.\\
Write $\phi_l$ the $\T^d_R$-periodic solution of the ergodicity problem 
$L_U\phi_{l}=|l-\nabla \chi_l|^2 - {^tl}D(U)l$
 with $\phi_l(0)=0$.
Write $F_l(x)=l.x-\chi_l(x)$ and $\psi_l(x)=F_l^2(x)-\phi_l(x)$, 
observe that since
$L_U F_l^2=|l-\nabla \chi_l|^2$ it follows that
\begin{equation}
L_U \psi_l={^tl}D(U)l
\end{equation}
The following inequality will be used to show that $\sum_{i=1}^d 
\psi_{e_i}$ behaves like $|x|^2$
\begin{equation}\label{AldIneq_Psi}
C_1 |x|^2 - C_2 (\|\chi_.\|_\infty^2+\|\phi_.\|_\infty) \leq 
\sum_{i=1}^d \psi_{e_i}(x) \leq
C_3 (|x|^2+\|\chi_.\|_\infty^2+\|\phi_.\|_\infty)
\end{equation}
Using theorem \ref{Stampa1} (theorem 5.4,  chapter 5 of \cite{St1}) to 
control $F_l$ and $\psi_l$ over one period (observing that $L_U F_l=0$, 
$L_U \psi_l=-1$) and using $\chi_l=l.x-F_l$ and $\phi_l=F_l^2-\psi_l$ 
one obtains easily that
$\|\phi_.\|_\infty \leq C_d e^{(9d+13)\Osc(U)} R^2$, combining this 
estimate with \eref{eqcontrchi9} one obtains 
\begin{equation}\label{aldcontrol_tau}
\|\chi_.\|_\infty^2+\|\phi_.\|_\infty \leq C_d e^{(9d+13)\Osc(U)} R^2
\end{equation}
Since $V$ has been assumed to be smooth (in a first step), 
 we can use Ito formula to obtain
\begin{equation}
\psi_l(y_t)=\int_{0}^{t} \nabla \psi_l(y_s) d\omega_s \;+\; 
{^tl}D(U)l\, t
\end{equation}
Now  ($e_i$ being an orthonormal basis of $\R^d$) write $M_{t}$ the 
local martingale
\begin{equation}
M_{t}=\sum_{i=1}^d \psi_{e_i}(y_t)-\operatorname{Trace}\big(D(U)\big)\, 
t
\end{equation}
Define
\[\tau'(0,r)=\inf \{t \geq 0 \;:\; |\sum_{i=1}^d \psi_{e_i}(y_t)|=r 
\}\]
According to the inequality \eref{AldIneq_Psi}
one has
\begin{equation}\label{aldtau1}
\tau'(0,C_1 r^2 - C_2 (\|\chi_.\|_\infty^2+\|\phi_.\|_\infty))
 \leq \tau(0,r) \leq \tau'(0,C_3 
(r^2+\|\chi_.\|_\infty^2+\|\phi_.\|_\infty))
\end{equation}
Since $M_{t \wedge \tau'(0,r)}$ is uniformly integrable (easy to prove 
by using the inequalities \eref{aldtau1})
one  obtains
\begin{equation}
\E[\tau'(0,r)]=\frac{r}{\operatorname{Trace}\big(D(U)\big)}
\end{equation}
Thus, by using the inequality \eref{aldcontrol_tau} and the 
Voigt-Reiss' inequality
$D(U) \geq e^{-2 \Osc(U)}$
one obtains the equation \eref{Aldtau_alphaeq3}.

\paragraph{}\label{paerxtihshsns7h53}
The proof of the weak stability result \eref{eq6g26hjsnm81} is based on 
the following obvious lemma that describes the monotony of Green 
functions as quadratic forms, i.e.
\begin{Lemma}\label{Tiwestarethprth1}
Let $\Omega$ be a smooth bounded open subset of $\R^d$. Assume that 
$M,Q$ are symmetric smooth coercive matrices on $\overline{\Omega}$. 
Assume $M \leq \lambda Q$ with $\lambda >0$, then for all $f \in 
C^0(\Omega)$, writing $G_Q$ the Green functions of $-\nabla Q \nabla$  with 
Dirichlet condition on $\partial \Omega$
\begin{equation}
\int_{\Omega}G_Q(x,y)f(y)f(x)\,dx\,dy \leq \lambda 
\int_{\Omega}G_M(x,y)f(y)f(x)\,dx\,dy
\end{equation}
\end{Lemma}
\begin{proof}
 Let $f \in C^0(\overline{\Omega})$. Write
$\psi_M$ , $\psi_Q$ the solutions of $-\nabla M \nabla \psi_M=f$ and 
$-\nabla Q \nabla \psi_Q=f$ with Dirichlet conditions on $\partial 
\Omega$.
Observe that $\psi_M$ and $\psi_Q$ are the unique minimizers of 
$I_M(h,f)$ and $I_Q(h,f)$ with
\begin{equation}
I_M(h,f)=\frac{1}{2} \int_{\Omega} {^t\nabla h} M \nabla h 
dx-\int_{\Omega}h(x) f(x) \,dx
\end{equation}
and  $I_M(\psi_M,f)=-\frac{1}{2}\int_{\Omega}\psi_M(x) f(x) \,dx$.
Observe that since $M\leq \lambda Q$
\begin{equation}\label{Tiwestarethpreq1}
\begin{split}
I_M(h,f)\leq \lambda I_Q(h,\frac{f}{\lambda})
\end{split}
\end{equation}
and the minimum of the right member in the equation 
\eref{Tiwestarethpreq1} is reached at $\psi_Q/\lambda$.
It follows that
$\int_{\Omega}\psi_Q(x) f(x) \leq \lambda \int_{\Omega}\psi_M(x) f(x)$,
which proves the lemma.
\end{proof}
Then, the equation \eref{eq6g26hjsnm81} follows directly from this 
lemma by choosing $Q=e^{-2(U+P)}$, $M=e^{-2U}$ and observing that 
$\E_x^{U}[\tau(0,r)]=2\int_{B(0,r)}G_{e^{-2U}I_d}(x,y)e^{-2U(y)}dy$.

\subsubsection{Anomalous heat kernel tail: theorem 
\ref{ddshdidcsbiii1}}
\paragraph{}\label{eqlkkecleecuu1}
From the pointwise anomaly of the hitting times of theorem 
\ref{IntSMAldcontrol_taudx2} one can deduce the anomalous heat kernel tail by 
adapting a strategy used by M.T. Barlow and R. Bass for the Sierpinski 
Carpet. This strategy is described in details in the proof of theorem 3.11 
of \cite{Ba98} and we will give only the main lines of its adaptation.\\
We will estimate $\P_x[\tau(x,r)<t]$ and use $\P_x[|y_t|>r]\leq 
\P_x[\tau(x,r)<t]$ to obtain theorem \ref{ddshdidcsbiii1}.\\
 Using the notations introduced in theorem \ref{IntSMAldcontrol_taudx2} 
and $M:=(d,K_0,K_\alpha,\alpha,\mu,\lambda_{\max})$, it will be shown 
in paragraph \ref{eqlkkecleecuu2}
 that for $r>C(M,\rho_{\max})$ one has
\begin{equation}\label{gsiasugsiiiuz1}
\P_x[\tau(x,r)\leq t] \leq 
\frac{t}{r^{2+\sigma(r)(1+\gamma)}C_{19}(M)}+1-C_{20}(M)r^{-2\gamma \sigma(r)}
\end{equation}
Now we will use lemma 3.14 of \cite{Ba98} given below (this is also 
lemma 1.1 of \cite{BB90}).
\begin{Lemma}\label{barlow_xi}
Let $\xi_1, \xi_2, \ldots,\xi_n$, V be non-negative r.v. such that 
$V\geq \sum_{i=1}^{n}\xi_i$.
Suppose that for some $p\in(0,1),\; a>0$ and $t>0$
\[\P\big(\xi_i \leq t|\sigma(\xi_1,\ldots,\xi_{i-1})\big) \leq p+at\]
Then
\[\ln \P(V \leq t) \leq 2 \big(\frac{ant}{p}\big)^{\frac{1}{2}}-n \ln 
\frac{1}{p}\]
\end{Lemma}
Let $n\geq 1$ and $g=\frac{r}{n}$.
Define the stopping times $S_i\; i\geq 0$ by $S_0=0$ and
$$S_{i+1}=\inf \{t\geq S_i \;:\;|y_t-y_{S_i}|\geq g\}$$
Write $\xi_i=S_i-S_{i-1}$ for $i\geq 1$
Let $\mathcal{F}_t$ be the filtration of $y_t$ and let 
$\mathcal{G}_i=\mathcal{F}_{S_i}$
Then it follows from the equation  \eref{gsiasugsiiiuz1} that for
\begin{equation}\label{gexvixeiwetx321}
r/n>C(M,\rho_{\max})
\end{equation}
\[\begin{split}
\P_x[\xi_{i+1}\leq t | \mathcal(G)_i]&=\P_{y_{S_i}}[\tau(y_{S_i},g)\leq 
t]\\
&\leq C_{21}(M) \frac{t}{g^{2+\sigma(r)(1+\gamma)}}+1-C_{20}(M) 
g^{-2\sigma(r)\gamma}
\end{split}\]
Since $|y_{S_i}-y_{S_{i+1}}|=g$ it follows that $\P_x$ a.s. 
$|x-y_{S_n}| \leq r$.
Thus
\[S_n=\sum_{i=1}^n \xi_i \leq \tau(x,r)\]
And by lemma \ref{barlow_xi} with
\[a=C_{21}(M)(\frac{n}{r})^{2+\sigma(r)(1+\gamma)}\quad p=1-C_{20}(M) 
(\frac{n}{r})^{2\sigma(r)\gamma}\]
One obtains
\begin{equation}\label{ujbbbbiqiw1}
\ln \P_x[\tau(x,r)\leq t] \leq 2 
\Big(\frac{n\,t\,C_{21}(\frac{n}{r})^{2+\sigma(r)(1+\gamma)}}{1-C_{20} (\frac{n}{r})^{2\sigma(r)\gamma}
}\Big)^{\frac{1}{2}}-n \ln \frac{1}{1-C_{20} 
(\frac{n}{r})^{2\sigma(r)\gamma}}
\end{equation}
Minimizing the right term in \eref{ujbbbbiqiw1} over $n$ under the 
constraint \eref{gexvixeiwetx321} and the assumptions 
\eref{gvisuziuzvvvv1}, $\rho_{\min}>C_{6,M}$, one obtains theorem \ref{ddshdidcsbiii1}.

\paragraph{}\label{eqlkkecleecuu2}
The equation \eref{gsiasugsiiiuz1} is an adaptation of lemma 3.16 of 
\cite{Ba98}.
Observe that
\[
\begin{split}
\E_x[\tau(x,r)] \leq &t+ \E_x[1(\tau(x,r)>t)\E_{y_t}[\tau(x,r)-t]]\\
\leq &t+ \P_x[1(\tau(x,r)>t)] \sup_{y\in B(x,r)}\E_{y}[\tau(x,r)]
\end{split}
\]
Using $\forall y\in B(x,r)$, $\P_y$ a.s. $\tau(x,r) \leq \tau(y,2r)$ it 
follows by theorem \ref{IntSMAldcontrol_taudx2} for 
$r>C(M,\rho_{\max})$
\[\begin{split}C_{33}(M) r^{2+\sigma(r)(1-\gamma)} \leq \E_x[\tau(x,r)] 
\leq t+\P_x[\tau(x,r)>t] 
C_{34}(M)r^{2+\sigma(r)(1+\gamma)}\end{split}\]
Which leads to \eref{gsiasugsiiiuz1}.

\paragraph*{Acknowledgments}
The authors would like to thank Alain-Sol Sznitman, Stefano Olla and 
Alano Ancona \cite{An99} for stimulating discussions.

\end{document}